%% file: paper.tex
\documentclass[11pt]{article}
\RequirePackage[OT1]{fontenc}
\RequirePackage{amsthm,amsmath,natbib}
\RequirePackage[colorlinks]{hyperref}
\RequirePackage{hypernat}
\usepackage{amsfonts}
\usepackage{amssymb}
\usepackage{amstext}
\usepackage{fullpage}
\bibpunct{(}{)}{;}{a}{,}{,}

\usepackage{times}
\usepackage{graphicx}
\usepackage{epsfig}

\input{./macros}

\def\supp{\mathop{\text{supp}\kern.2ex}}
\def\argmin{\mathop{\text{arg\,min}\kern.2ex}}
\let\hat\widehat
\let\tilde\widetilde

\def\Z{\widetilde X}

\numberwithin{equation}{section}
\theoremstyle{plain}

\newtheorem{thm}{Theorem}[section]
\newtheorem{assu}{Assumption}[section]
\newtheorem{lmma}[thm]{Lemma}
\newtheorem{propn}[thm]{Proposition}
\newtheorem{define}[thm]{Definition}
\newtheorem{remrk}[thm]{Remark}
\newtheorem{clm}[thm]{Claim}
\newtheorem{coroll}[thm]{Corollary}
\newenvironment{proposition}{\begin{propn}\hskip-6pt{\bf .}\enspace \sl}{\end{propn}}
\newenvironment{theorem}{\begin{thm}\hskip-6pt{\bf .}\enspace \sl}{\end{thm}}
\newenvironment{assumption}{\begin{assu}\hskip-6pt{\bf .}\enspace \sl}{\end{assu}}
\newenvironment{lemma}{\begin{lmma}\hskip-6pt{\bf .}\enspace \sl}{\end{lmma}}

\newenvironment{definition}{\begin{define}\hskip-6pt{\bf .}\enspace \sl}{\end{define}}
\newenvironment{corollary}{\begin{coroll}\hskip-6pt{\bf .}\enspace \sl}{\end{coroll}}
\newenvironment{remark}{\begin{remrk}\hskip-6pt{\bf .}\enspace \rm}{\end{remrk}}
\def\v#1{{{\mbox{{\mbf #1}}}\kern.2ex}}

\newenvironment{subeqnarray}
{\begin{subequations}\begin{eqnarray}}%
{\end{eqnarray}\end{subequations}\hskip-4.0pt}%
\begin{document}
\mbox{\ }
\vskip.2in
\centerline{\LARGE\bf Adaptive Lasso for High Dimensional Regression and}
\vskip .05in
\centerline{\LARGE\bf Gaussian Graphical Modeling \footnote{Research supported by SNF 20PA21-120050/1.}} 
\vskip .2in
\begin{center}
\begin{tabular}{c}
{\large Shuheng Zhou\;\; Sara van de Geer\;\; Peter B\"{u}hlmann} \\[15pt]
Seminar f\"{u}r Statistik \\
ETH Z\"{u}rich \\
CH-8092 Z\"{u}rich, Switzerland \\[20pt]
\end{tabular}
\end{center}
\begin{center}
March 13th, 2009
\end{center}

\begin{abstract}
We show that the two-stage adaptive Lasso procedure \citep{Zou06} is
consistent for 
high-dimensional model selection in linear and Gaussian graphical
models. Our conditions for consistency cover more general situations 
than those accomplished in previous work: we prove that 
restricted eigenvalue conditions \citep{BRT08} are also
sufficient for sparse structure estimation.  
\end{abstract}




\section{Introduction}
\label{sec:introduction}

The problem of inferring the sparsity pattern, i.e. model selection, in
high-dimensional problems has recently gained a lot of attention. One
important stream of research, which we also adopt here, requires
computational feasibility and provable statistical properties of estimation
methods or algorithms. 
Regularization with $\ell_1$-type penalization has become extremely
popular for model selection in high-dimensional scenarios. The methods are
easy to use, due to recent progress in convex optimization
\citep{mevdgpb08}, \citep{fht08}, and they are  
asymptotically consistent or oracle optimal when requiring some conditions,
e.g. on the design 
matrix in a linear model or among the variables in a graphical
model \citep{GR04,MB06,vandeG08}, \citep{BRT08}. However, these conditions, 
referred to as coherence or 
compatibility conditions, are often very restrictive. The restrictions are
due to severe bias problems with $\ell_1$-penalization, i.e. shrinking also
the estimates which correspond to true signal variables,
see also \cite{Zou06}, \cite{ME07}. 

Regularization with the $\ell_q$-norm with $q <1$ would mitigate
some of the bias problems but become computationally infeasible as the
penalty is non-convex. As an interesting alternative, one can consider
multi-step procedures where each of the steps involves a convex
optimization only. A prime example is the adaptive Lasso \citep{Zou06}
which is a 
two-step algorithm and whose repeated application corresponds in some
``loose'' sense to a non-convex penalization scheme \citep{zouli2008}. We are
analyzing in this paper this adaptive Lasso procedure for variable
selection in linear models as well as for Gaussian graphical modeling. Both
frameworks are related to each other and for both of them, we derive
results for model selection under rather weak conditions. In
particular, our results imply that the adaptive Lasso can recover the true
underlying model in situations where plain $\ell_1$-regularization fails
(assuming restricted eigenvalue conditions).  

\subsection{Variable selection in linear models}

Consider the linear model 
\beq
\label{eq::linear-model}
Y = X \beta + \epsilon,
\eeq
where $X$ is an $n \times p$ design matrix, $Y$ is an $n \times 1$ vector
of noisy observations and $\epsilon$ being the noise term. The design
matrix is treated as either fixed or random. 
We assume throughout this paper that $p \ge n$ (i.e. high-dimensional) and 
$\epsilon \sim N(0, \sigma_{\epsilon}^2 I_n)$. 

The sparse object to recover is the unknown parameter $\beta \in \R^p$. We
assume that 
it has a relatively small number $s$ of nonzero coefficients:
$S := \supp(\beta)$ \break $ = \left\{j \,:\, \beta_j \neq 0\right\}$ and  
$s = \abs{\supp(\beta)}$.
Let $\beta_{\rm min} := \min_{j \in S} |\beta_{j}|$.
Inferring the sparsity pattern, i.e. variable selection, 
refers to the task of correctly estimating the support set
$\supp(\beta)$ based on noisy observations from (\ref{eq::linear-model}). 
In particular, given some estimator $\hat\beta$, recovery of the relevant
variables is understood to be 
\begin{eqnarray}\label{eq:signcons}
{\rm supp}(\hat\beta) = {\rm supp}(\beta)\ \mbox{with high probability}. 
\end{eqnarray}

Regularized estimation with the $\ell_1$-norm penalty, also 
known as the Lasso \citep{Tib96}, refers to the following convex
optimization problem:   
\begin{eqnarray}
\label{eq::origin}
\hat \beta = \arg\min_{\beta} \frac{1}{2n}\|Y-X\beta\|_2^2 + \lambda_n \|\beta\|_1,
\end{eqnarray}
where the scaling factor $1/(2n)$ is chosen by convenience and $\lambda_n \ge
0$ is a penalization parameter. It is an attractive and computationally
tractable method with provable good statistical properties, even if $p$ is
much larger than $n$, for prediction
\citep{GR04}, for estimation in terms of the $\ell_1$- or $\ell_2$-loss
\citep{vandeG08,MY08,BRT08} and for variable selection
\citep{MB06,ZY07,Wai08}. For the
specific problem of variable selection, it is known that the so-called
``neighborhood stability condition'' for the design matrix \citep{MB06},
which has been re-formulated in a nicer form as the ``irrepresentable condition''
\citep{ZY07}, is necessary and sufficient for consistent variable selection
in the sense of (\ref{eq:signcons}). Moreover, as this condition is
restrictive, its necessity implies that the Lasso only works in a
rather restricted range of problems, excluding cases where the design
exhibits too strong (empirical) correlations.
A key motivation of our work is to continue the exploration of a 
computationally tractable algorithm for variable selection, 
while aiming to relax the stringent conditions that are imposed on the 
design matrix $X$. 

Towards these goals, we analyze the adaptive Lasso procedure, see
(\ref{eq::weighted-lasso}) below, for 
variable selection in the high-dimensional setting. This method was originally
proposed by~\cite{Zou06} and he analyzed the case when $p$ is
fixed. Further progress of analyzing the adaptive Lasso in the
high-dimensional scenario has been achieved 
by \cite{HMZ08}. A more complete understanding of its power, when applied to 
the high dimensional setting where $p \gg n$ is still lacking. We prove in
this paper that variable selection with the adaptive Lasso is possible
under rather general incoherence conditions on the design. We do not
require more stringent conditions on the design $X$ than \cite{BRT08} who
give the currently weakest conditions for convergence of the Lasso in terms of 
$\|\hat{\beta} - \beta\|_1$ and $\|\hat{\beta} - \beta\|_2$. We show that
for an initial estimator $\beta_{\init}$ in the two-stage adaptive Lasso
procedure with a sufficiently reasonable behavior of $\|\beta_{\init} -
\beta\|_{\infty}$, model selection is possible assuming only a lower bound
on the smallest eigenvalue of $X_S^T X_S/n$, where $X_S$ denotes the
submatrix of $X$ whose columns are indexed by $S$, and some restrictions on
$\beta_{\min}$ and the sparsity level $s$. 
Thus, variable selection is possible under rather general design conditions 
by the two-stage adaptive Lasso, and it is necessary to move away from
plain $\ell_1$-regularization, see \cite{MB06}, \cite{ZY07}.  

\subsection{Covariance selection in Gaussian graphical models}\label{subsec.covsel}

Covariance selection in a Gaussian graphical model refers to the problem of
inferring conditional independencies between a set of jointly Gaussian
random variables
\begin{eqnarray}\label{eq:multGauss}
X_1,\ldots X_p \sim N(0,\Sigma)
\end{eqnarray}
(the restriction to mean 0 is without loss of generality). 
These variables $X_1,\ldots ,X_p$ correspond to nodes in a graph, labeled
by $\{1,\ldots ,p\}$, and a
Gaussian conditional independence graph is then defined as follows:
\begin{eqnarray*}
\mbox{there is an undirected edge between node $i$ and $j$}\
\Leftrightarrow\ \Sigma^{-1}_{ij} \neq 0.
\end{eqnarray*}
The definition of an edge is equivalent to requiring that $X_i$ and $X_j$ are
conditionally dependent given all remaining variables $\{X_k;\ k \neq
i,j\}$. For details cf. \cite{laur96}. Estimation of the edge set is thus
equivalent to finding the zeroes in the concentration matrix
$\Sigma^{-1}$. 

In the high-dimensional scenario with $p \ge n$, where $n$
denotes the sample size of i.i.d. copies from (\ref{eq:multGauss}),
$\ell_1$-type regularization has been analyzed. \break
\cite{MB06} prove that it
is possible to consistently infer the edge set by considering many variable
selection problems in high-dimensional Gaussian regressions, again
requiring a global neighborhood stability or irrepresentable condition
which puts some restrictions on the covariance matrix $\Sigma$. Later, the
GLasso penalization has been proposed \citep{FHT07,BGD07} which is a sparse
estimator for 
$\Sigma^{-1}$ using an $\ell_1$-penalty on the non-diagonal elements of
$\Sigma^{-1}$ in the multivariate Gaussian log-likelihood. 
~\cite{RWRY08} recently obtained results for consistent covariance 
selection ((i.e. inferring the edge set) using the GLasso by
imposing mutual incoherence conditions (analogous to the neighborhood
stability condition) on the Fisher information matrix 
(of size $p^2 \times p^2$) of the model, which is  an edge-based counterpart 
of $\Sigma$.

We focus here on generalizing conditions for the pursuit via many regressions:
we prove in this paper a result for inferring the edge set in a Gaussian 
graphical model, under a rather general condition on $\Sigma$ closely 
related to the restricted eigenvalue assumptions in \cite{BRT08}
by analyzing the pursuit of many regressions with the adaptive Lasso. 
We conjecture that the set of conditions which we are imposing are 
more general than what~\cite{RWRY08} require when using the GLasso, 
although this is a point that needs to be thoroughly studied as we discuss 
further in Section~\ref{sec.discussion}.
We also suspect that the GLasso approach is intrinsically more limited, 
in terms of restrictions for the covariance matrix $\Sigma$ than the approach 
from \cite{MB06} via considering many regressions. 
This has been recognized by \cite{meinshausen05} and also studied 
by~\cite{RWRY08} on specific graphical models.
On the other hand, for well-behaved problems, GLasso might
have an advantage because it exploits the positive definiteness of $\Sigma$
and $\Sigma^{-1}$.

\subsection{Related work}
\label{sec:background}
Recently,~\cite{HMZ08} studied the adaptive Lasso
estimators in sparse, high-dimensional linear regression models for a 
fixed design. Under a rather strong mutual incoherence condition between 
every pair of relevant and irrelevant covariates and assuming other regularity
conditions, they prove that the adaptive Lasso recovers the correct  
model and has an oracle property.
While they have derived the same incoherence condition as one (among
others) of ours in~\eqref{eq:eta} in order for the second stage 
weighted Lasso procedure to achieve model selection 
consistency, they achieve it by an initial estimator assuming some strong
mutual incoherence condition which bounds the pairwise correlations of the
columns of the design. 
This is a much stronger condition than the 
restricted eigenvalue assumptions that we make, see~\cite{BRT08}.

\cite{MY08} examined the variable selection property of the Lasso followed
by a thresholding procedure. Under a relaxed ``incoherence design'' assumption, \cite{MY08}
show that the estimator is still consistent in the $\ell_2$-norm sense for 
fixed designs, and furthermore, it is possible to do
hard-thresholding on the ordinary Lasso estimator to achieve variable
selection consistency. However the choice of the threshold parameter
depends on the the unknown value $\beta_{\min}$ and the sparsity $s$ of 
$\beta$. It is not clear how one can choose
such a threshold parameter without knowing $\beta_{\min}$ or $s$. 
A more general framework for multi-stage variable selection 
was studied by~\cite{WR08} for various methods and conditions. Their
approach controls the probability of false positives (i.e. type I error) but
pays a price in terms of false negatives (i.e. type II error) in comparison
to the adaptive Lasso \citep{WR08}.

Finally, our focus is rather different from that of~\cite{Wai08,Wai07},
where the goal was to analyze the least 
amount of samples that one needs in order to recover a sparse signal via a 
random or a fixed measurement ensemble that satisfies strong incoherence 
conditions. It is an open problem to establish a lower bound on the sample
size, given $p$, $s$ and $\beta_{\min}$, to recover the model with the
adaptive Lasso, assuming restricted eigenvalue assumptions only.  

\subsection{Organization of the paper}
In Section~\ref{sec:general}
we define the two-step adaptive Lasso procedure for linear regression and
describe our main result: general model selection
properties of the second stage weighted procedure for variable
selection. Here, the initial estimator $\beta_{\rm init}$
can be general, and we assume
a bound for $ \| \beta_{\rm init} - \beta \|_{\infty}$.
Section ~\ref{sec:summary} presents the restricted eigenvalue conditions
we need for deriving bounds for 
$\|  \beta_{\init} - \beta\|_{\infty} $ with the standard Lasso  as initial estimator $\beta_{\init}$.
In Sections~\ref{adapfix},  \ref{adaprandom} and~\ref{sec:graph}, 
we summarize conditions and results, with the standard Lasso as
initial estimator,  for linear regression with fixed design,
linear regression with random design, and 
Gaussian graphical modeling, respectively. These results are consequences
of our general result in Section~\ref{sec:general}.
Section~\ref{sec:sparsistence} presents a model selection lemma
for the weighted Lasso with general weights.
The remainder of the paper contains the proofs. 

\section{The adaptive Lasso estimator and its general properties} 
\label{sec:general}
Consider the linear model in (\ref{eq::linear-model}). We distinguish
later between fixed and random design.
\subsection{The two-stage adaptive Lasso procedure}
\label{sec:adap-procedure}

The adaptive Lasso is the Lasso estimator with a re-weighted penalty function,
see (\ref{eq::weighted-lasso}) below. The weights are estimated from an
initial estimator $\beta_{\init}$: 
\beq
\label{eq::adaptive-weights}
w_j := \max \{\inv{|\beta_{j,{\rm init}}|}, 1\} .
\eeq 
We note that the original proposal of \cite{Zou06} uses $w_j =
1/|\beta_{j,{\rm init}}|^{\gamma}$ for some $\gamma > 0$ with $\gamma = 1$
the most common choice.
The adaptive Lasso is now defined by a second-stage weighted Lasso: 
\beq
\label{eq::weighted-lasso}
\hat \beta = \arg \min_{\beta}  
\inv{2n} \| Y- X\beta \|^2_2 + \lambda_n \sum_{j=1}^p w_j | \beta_j | . 
\eeq

\subsection{Variable selection with the adaptive Lasso estimator}
\label{sec:intro-weighted}
Correct variable selection with the adaptive Lasso requires some conditions
for the design. We first make some assumptions related to the
design matrix. For a symmetric matrix $A$, let $\Lambda_{\min}(A)$ denote 
the  smallest eigenvalue of $A$.

For a fixed design matrix $X$, we define
\beq
\label{eq::eigen-admissible}
\Lambda_{\min}(s) := \min_{\stackrel{J_0 \subseteq \{1, \ldots, p\}}{|J_0|
    \leq s}} 
\min_{\stackrel{\gamma \not=0}{\gamma_{J_0}^c = \; 0}}
\; \;  \frac{\twonorm{X \gamma}^2}{n\twonorm{\gamma_{J_0}}^2}.
\eeq
We assume throughout this paper that $\Lambda_{\min}(s) > 0$ . As a
consequence of this definition we have,
\begin{eqnarray}
\label{eq::eigen-cond}
\Lambda_{\min}  \left(\frac{X_S^T X_S}{n}\right) \geq \Lambda_{\min}(s) > 0.
\end{eqnarray}
Furthermore, we assume for fixed design that the $\ell_2$-norm of each
column of $X$ is upper bounded by $c_0 \sqrt{n}$ for some constant $c_0
>0$. We then consider the set 
\beq
\label{eq::low-noise}
{\T} := 
\biggl \{\norm{\frac{X^T \e}{n}}_{\infty}
\le c_0 \sigma_{\e}\sqrt{\frac{6\log p}{n}} \biggr \}.
\eeq
The set $\T$ has large probability, as described below in (\ref{eq:set-T}). 

For a random design matrix $X$ we assume:
\begin{eqnarray}
\label{eq::rand-des}
X\ \mbox{has i.i.d. rows}\ \sim N(0, \Sigma),
\end{eqnarray}
where we assume without loss of
generality that the mean is zero and $\Sigma_{jj} = 1, \forall j = 1, \ldots,
p$. We then define 
\beq
\label{eq::eigen-admissible-random}
\Lambda_{\min}(s) :=
\frac{16}{17} \min_{\stackrel{J_0 \subseteq \{1, \ldots, p\}}{|J_0| \leq s}}
\min_{\stackrel{\gamma \not=0}{\gamma_{J_0}^c = \; 0}}
\; \;  \frac{\gamma^T \Sigma \gamma}{\twonorm{\gamma_{J_0}}^2}.
\eeq
As for fixed design, we assume that $\Lambda_{\min}(s) >0$ with large
probability. 
The factor $16/17$ allows us to use the same notation
$\Lambda_{\min}(s)$ for both fixed and random design. Let $\Sigma_{SS}$ be
the sub-matrix with rows and columns both indexed by  
the active set $S$. It then holds that 
\begin{eqnarray}\label{eq::eigen-cond-random}
\Lambda_{\min} (\Sigma_{SS} ) & \geq  & \frac{17 \Lambda_{\min}(s)}{16} >0.
\end{eqnarray}
Then, a random design $X$ as in (\ref{eq::rand-des}) behaves nicely, with high
probability. To be more precise, denote by $\Delta = \frac{X^T X}{n} -
\Sigma$, and consider 
\begin{eqnarray}
\label{eq::good-random-design}
\X := 
\left\{\max_{j, k} |\Delta_{jk}| < C_2 \sqrt{\frac{\log p}{n}} \right\},
\end{eqnarray} 
for some constant $C_2> 4 \sqrt{5/3}$.
Throughout this paper, we assume for random design that 
$p < e^{n/4C_2^2}$, i.e. $C_2 \sqrt{\frac{\log p}{n}} < 1/2$, such that
$\X$ holds with probability at least  
$1 - \frac{1}{p^2}$ (cf.~\eqref{eq::X-bound} and 
Lemma~\ref{lemma:gaussian-covars}). We note that this implies that on $\X$,
\begin{eqnarray}
\label{eq::good-random-design-diag}
\forall j = 1, \ldots, p, \; \; \; \twonorm{X_j}^2 \leq \frac{3n}{2}.
\end{eqnarray}
The set $\T$ in (\ref{eq::low-noise}), intersected with
$\X$, is also relevant for random design: the constant $c_0$ equals
$\sqrt{3/2}$, following~\eqref{eq::good-random-design-diag}.

For both, fixed and random design, we consider the quantity
\begin{eqnarray}
\label{eq::r_n_bound}
r_n(S) := \norm{X^T_{S^c} X_S (X_S^T X_S)^{-1}}_{\infty}, 
\end{eqnarray}
where $\|A\|_{\infty} = \max_{1 \le i \le k}
\sum_{j=1}^m |A_{ij}|$ for a $k \times m$ matrix $A$. 
The properties of the adaptive Lasso procedure depend on (an upper bound
of) $r_n(S)$.  

Finally, we denote by 
$$\delta := \beta_{\init} - \beta$$ 
the difference
between the initial estimate and the true parameter value. 

\begin{theorem}
\label{thm:first}
Consider the adaptive Lasso estimator in a linear model as in
(\ref{eq::linear-model}) with design $X$, where $n \le p$, and for 
fixed design: the $\ell_2$-norm of each
column of $X$ is upper bounded by $c_0 \sqrt{n}$ for some constant $c_0 >0$.

Assume the upper bound
$\tilde{r_n} \ge r_n(S)$ which we  
require to hold only on $\X$ in
case of a random design. Furthermore, assume on  
$\T$ for a fixed design and on $\X \cap \T$ for a random design, some
upper bounds on $\delta$ as follows:  
$1 > \tilde{\delta}_S \ge \|\delta_S\|_{\infty}$ and
$1 > \tilde{\delta}_{\Sc} \ge \|\delta_{\Sc}\|_{\infty}$. 
Suppose that on
$\T$ for a fixed design and on $\X \cap \T$  
for a random design:\\
for some $1> \eta  > 0$ and some  
constant $M \geq \frac{4}{\eta}$, $\lambda_n$ is chosen from the range
\begin{equation}
\label{eq::penalty-choice}
M c_0 \sigma \tilde{\delta_{\Sc}}
\sqrt{\frac{2 \log(p-s)}{n}} \geq \lambda_n 
\geq \frac{4 c_0 \sigma \tilde{\delta_{\Sc}}}{\eta} 
\sqrt{\frac{2\log(p-s)}{n}}.
\end{equation}
Furthermore, assume:
\begin{eqnarray}
\label{eq::sparsity-condition-first}
\ \ \ \tilde{r_n} \leq 
\frac{1-\eta}{\tilde{\delta}_{\Sc}},
\end{eqnarray}
and for $C_1 = \max\left\{\frac{2\tilde{r_n}}{1-\eta},
  \frac{M}{\sqrt{3}}\right\}$
\begin{eqnarray}
\label{eq::tragic-0}
\ \ \ \ \ \ \ \ \ \ \beta_{\min}  >
\max\left\{2 \tilde{\delta}_{S},
\frac{2\lambda_n \sqrt{s}}{\Lambda_{\min}(s)},
\frac{4 c_0 \sigma}{\Lambda_{\min}(s)} \sqrt{\frac{6 s \log p}{n}},
C_1 \tilde{\delta}_{\Sc}\right\}.
\end{eqnarray}
Then, with probability $1 - \prob{\T^c} - 1/p^2$ for a fixed design or $1 -
\prob{(\X \cap \T)^c} - 1/p^2$ for a random design respectively, the
optimal solution $\hat{\beta}$ to (\ref{eq::weighted-lasso}) satisfies 
${\rm supp}(\hat{\beta}) = {\rm supp}(\beta)$.   
\end{theorem}
A proof is given in Section~\ref{finalproof}. We furthermore argue below
that the sets $\T$ and $\X$ (and hence also $\T \cap \X$) have large
probability.  
 
\begin{remark}\label{rem:unique}
In general, there are multiple solutions of the adaptive Lasso in
(\ref{eq::weighted-lasso}). However, with high probability, the solution of
(\ref{eq::weighted-lasso}) is unique. This follows from \cite{Wai08} and we
present more details in Section \ref{subsec.uniqueness}.  
\end{remark}

\begin{remark}
The last term on the right hand side in~\eqref{eq::tragic-0} usually
dominates all others (under the assumptions we make for the theorem): the
order of magnitude is typically $O(\sqrt{s \log(p) /n})$. Furthermore, for
a fixed design, we emphasize that 
$\tilde{r}_n$, $\tilde{\delta}_{S}$ and
$\tilde{\delta}_{\Sc}$ are only required  to hold on the 
set $\T$. Similarly, for a random design, we only require some upper bounds
to hold on the set $\T \cap \X$. 
\end{remark}

\begin{remark}
We note that Theorem \ref{thm:first} suggests that we can use any initial
estimator 
that yields a nice bound on 
$\norm{\delta}_{\infty} = \norm{\beta_{\init} - \beta}_{\infty}$. We
consider the Lasso as initial estimator in Sections \ref{adapfix} and
\ref{adaprandom}. The
Dantzig selector \citep{CT07} could be an alternative having similar 
properties as the Lasso under the restricted eigenvalue 
assumptions \citep{BRT08}. 
\end{remark}

\begin{lemma}\label{lemma:twonorm-random}
For a fixed design,  we have 
\begin{eqnarray}\label{eq:set-T}
\prob{\T} \ge 1 - 1/p^2.
\end{eqnarray}
Moreover, 
for a random design $X$ as in (\ref{eq::rand-des}) with 
$\Sigma_{jj} = 1, \forall j \in \{1, \ldots, p\}$,
and for $p < e^{n/4C_2^2}$, where $C_2 > 4 \sqrt{5/3}$,
we have
\begin{eqnarray}
\label{eq::X-bound}
\prob{\X} \ge 1 - 1/p^2.
\end{eqnarray}
Hence, for a random design,
\begin{eqnarray*}
\prob{\X \cap \T} \ge 1 - 2/p^2.
\end{eqnarray*}
\end{lemma}
A proof is given in Section \ref{subsec:lemmatwonorm} (Lemmas
\ref{lemma:gaussian-noise} and  \ref{lemma:gaussian-covars}).

\section{Restricted eigenvalue conditions}\label{sec:summary}

We are analyzing in later sections the properties of the adaptive Lasso when
using the standard Lasso as initial estimator: 
\beq
\label{eq::initial-estimator}
\beta_{\init} :=  \arg \min_{\beta}  
\inv{2n} \| Y- X\beta \|^2_2 + \lambda_{\init} \sum_{j=1}^p | \beta_j |,
\eeq
where for some constant $B$ and $c_0$ to be specified,
\beq
\label{eq::init-lower-bound}
\lambda_{\init} = B c_0 \sigma_{\e} \sqrt{\frac{\log p}{n}}.
\eeq
As usual, in order to be a sensible procedure, we assume that the different
variables (columns in $X$) are on the same scale.
In view of Theorem \ref{thm:first},
we need to establish bounds for $\delta = \beta_{\init} -\beta$, where
$\beta_{\init}$ is defined in (\ref{eq::initial-estimator}). 

To derive such bounds for $\delta$, we build upon recent work
by~\cite{BRT08} under the ``restricted eigenvalue'' assumptions 
formalized therein, which are weaker than those in~\cite{CT07,MY08} for 
deriving $\ell_p$ bounds on $\delta$, where $p=1, 2$, for the Dantzig 
selector and the Lasso respectively. Similar conditions have been used 
by~\cite{Kol07} and~\cite{vandeG07}. 

\subsection{Restricted eigenvalue assumption for fixed design}
To introduce the first assumption, we need some more notation.
For integers $s, m$ such that $1 \leq s \leq p/2$ and $m \geq s, s + m \leq p$,
a vector $\delta \in \R^p$ and a set of indices $J_0 \subseteq \{1, \ldots, p\}$
with $|J_0| \leq s$, denoted by $J_m$ the subset of $\{1, \ldots, p\}$
corresponding to the $m$ largest in absolute value coordinates of $\delta$
outside of $J_0$ and defined $J_{0m} \stackrel{\triangle}{=} J_0 \cup J_m$.
\begin{assumption}\textnormal{{\bf Restricted eigenvalue assumption} $RE(s,
    m, k_0, X)$~\citep{BRT08}}. 
\label{def:BRT-cond-two}
Consider a fixed design. For some integer $1\leq s \leq p/2$, $m \geq s, s
+ m \leq p$, and 
a positive number $k_0$, the following condition holds:
\beq
\label{eq::admissible-two}
\inv{K(s, m, k_0, X)} := 
\min_{\stackrel{J_0 \subseteq \{1, \ldots, p\},}{|J_0| \leq s}}
\min_{\stackrel{\gamma \not=0,}
{\norm{\gamma_{J_0^c}}_1 \leq k_0 \norm{\gamma_{J_0}}_1}}
\; \;  \frac{\norm{X \gamma}_2}{\sqrt{n}\norm{\gamma_{J_{0m}}}_2} > 0.
\eeq
\end{assumption}
We often restrict ourselves to the case with $k_0 = 3$. 
Apparently, $RE(s, m, k_0, X)$ implies that
$RE(s, k_0, X)$ as in Definition~\ref{def:BRT-cond} below holds with 
$K(s, k_0, X) \leq K(s, m, k_0, X)$ for the same $X$.

\begin{definition}
\textnormal{{\bf Restricted eigenvalue definition} $RE(s, k_0,
  X)$~\citep{BRT08}}.  
\label{def:BRT-cond}
Consider a fixed design. For some integer $1\leq s \leq p$ and a positive
number $k_0$, the following condition holds:
\beq
\label{eq::admissible}
\inv{K(s, k_0, X)} :=
\min_{\stackrel{J_0 \subseteq \{1, \ldots, p\},}{|J_0| \leq s}}
\min_{\stackrel{\gamma \not=0,}{\norm{\gamma_{J_0^c}}_1 \leq k_0 
\norm{\gamma_{J_0}}_1}}
\; \;  \frac{\norm{X \gamma}_2}{\sqrt{n}\norm{\gamma_{J_0}}_2} > 0.
\eeq
\end{definition}
We note that variable selection with the adaptive Lasso is possible
under this weaker form of restricted eigenvalues, though with stronger
conditions on the sparsity $s$ and $\beta_{\min}$. 
We omit such results in this paper due to the 
lack of space.

By an argument in~\cite{BRT08}, it is known that  if $RE(s, k_0, X)$ is 
satisfied with $k_0 \geq 1$, then the square submatrices of size $\leq 2s$ 
of $X^T X/n$ are necessarily positive definite. 
In fact, it is clear that 
in~\eqref{eq::admissible}, the set of admissible $\gamma$ is a superset 
of that in~\eqref{eq::eigen-admissible}. Hence we have the following:

\begin{proposition}
\label{prop:eigen-K-bound}
Suppose Assumption $RE(s, k_0, X)$ holds for $1 \leq s \leq p/2$ and 
some $k_0>0$. Then $\Lambda_{\min}(s) \geq \inv{ K^2(s, k_0, X)} > 0$ for 
$\Lambda_{\min}(s)$ as defined in~\eqref{eq::eigen-admissible}.
\end{proposition}
Note that the quantity $\Lambda_{\min}(s)$ also appears in Theorem
\ref{thm:first} and hence when applying it, we make use of Proposition
\ref{prop:eigen-K-bound}.  

\subsection{Restricted orthogonality assumption for fixed design}
We also present results under a stronger design condition which covers
cases where the sparsity $s$ is allowed to be larger than in Corollary
\ref{thm:adaptive-recovery-two} under Assumption
\ref{def:BRT-cond-two}, see also Corollary
\ref{thm:adaptive-recovery-three}.   
We define the $(s, s')$-restricted orthogonality constant~\citep{CT07}
$\theta_{s, s'}$ for $s + s' \leq p$, which is the smallest quantity such that
\beq
\label{label:correlation-coefficient}
\abs{\frac{\ip{X_T c, X_{T'} {c'}}}{n}} 
\leq \theta_{s, s'} \twonorm{c} \twonorm{c'}
\eeq
holds for all disjoint sets $T,T' \subseteq \{1, \ldots, p\}$ 
of cardinality $|T| \leq s$ and $|T'| < s'$. 
\begin{assumption}
\textnormal{\bf{Restricted orthogonality assumption}}. 
\label{def:BRT-cond-three}
Consider a fixed design. For some integer $1\leq s \leq p/2$, $m \geq s, s
+ m \leq p$, and 
a positive number $k_0$, the condition $RE(s, s, k_0, X)$
holds. Furthermore, the following condition holds:
\begin{eqnarray}
\label{eq::admissible-three}
& &\Lambda_{\min}(s) >
16 k_0 K^2(s, m, k_0, X) \lambda_{\init} s \theta_{1, s},\\
& &\label{eq::linear-sparsity}
s < \frac{n}{96 c_0^2  \sigma^2 K^2(s, s, k_0, X) \log p},
\end{eqnarray}
where $k_0 \le 3$. 
\end{assumption}
%
With such a restriction on the sparsity, we note that
\eqref{eq::admissible-three} is a weaker condition than
Assumption $3$ in~\cite{BRT08}. 
We assume that ~\eqref{eq::admissible-three} holds with a constant that is
smaller than $2 k_0$ as in Assumption $3$ of~\citep{BRT08}, which by 
itself is a sufficient condition to derive
Assumption~\ref{def:BRT-cond}. 

We refer to~\cite{BRT08} for more detailed discussions about these
assumptions which are weaker than those in~\cite{CT07,MY08} and arguably 
less restrictive than those in~\cite{MB06},\\
~\cite{ZY07} or \cite{Wai08}.

\section{The adaptive Lasso with fixed design}\label{adapfix}
We first show that the restricted eigenvalue condition ensures to derive 
upper bounds on the $\ell_{\infty}$-norms of
$\delta := \beta_{\init} - \beta$, 

\begin{lemma}
\label{lemma:weights-bound}
Suppose that condition $RE(s, 3, X)$ holds for a fixed design and suppose that
\begin{eqnarray}
\label{eq::weight-bounds-conditions}
\beta_{\rm min} & \ge & 8 K^2(s, 3, X) \lambda_{\init} \sqrt{s},
\end{eqnarray}
for $\lambda_{\init}$ that satisfies~\eqref{eq::init-lower-bound}.
Then, the initial estimator~\eqref{eq::initial-estimator} in model 
~\eqref{eq::linear-model} guarantees that 
on the set $\T$ as in~\eqref{eq::low-noise},
\begin{subeqnarray}
\label{eq::e1}
\norm{\delta_{S}}_{\infty} & \leq & 
4 K^2(s, 3, X) \lambda_{\init} \sqrt{s}, \text{ and }\\ 
\label{eq::weight-bounds-min}
\norm{\delta_{\Sc}}_{\infty} & \leq &  3 K^2(s, 3, X) \lambda_{\init} s
\end{subeqnarray}
Suppose that Assumption $RE(s, s, 3, X)$
and~\eqref{eq::weight-bounds-conditions} hold. Then
on the set $\T$ as in~\eqref{eq::low-noise}, 
~\eqref{eq::e1} holds, 
while~\eqref{eq::weight-bounds-min} is replaced by
\begin{eqnarray}
\label{eq::weight-bounds-min-two}
\norm{\delta_{\Sc}}_{\infty} \leq  16 K^2(s, s, 3, X) \lambda_{\init} \sqrt{s}.
\end{eqnarray}
\end{lemma}
A proof is given in Subsection \ref{sec:initial}.
Lemma \ref{lemma:weights-bound} leads to the upper bounds
$\tilde{\delta}_S = 4 K^2(s, 3, X) \lambda_{\init} \sqrt{s}$ and
$\tilde{\delta}_{\Sc} = 16 K^2(s, s, 3, X) \lambda_{\init} \sqrt{s}$. 
When using 
these bounds in Theorem \ref{thm:first}, we see that the range for the
regularization parameter in \ref{eq::penalty-choice} depends on the unknown
sparsity $s$. This unpleasant situation can be improved by estimating $s$ 
using a thresholding procedure as follows. 

\begin{lemma}\textnormal{{\bf Thresholding procedure.}}
\label{thres} 
Let the assumptions of Lemma \ref{lemma:weights-bound} hold.
Consider the set $\bar{S}$ that includes all $\beta_{j, \init}$ for $j \in
\{1, \ldots,  p\}$,  whose absolute values are larger than $4
\lambda_{\init}$. Let $\bar{s} := |\bar{S}|$ be an estimate 
which is in the same order as the true sparsity $s$. More specifically,
we have, on the set $\T$ in~\eqref{eq::low-noise},   
\begin{equation}
\label{eq::superset}
S \subseteq \bar{S} \text{ and } \; \; 
s \leq |\bar{S}| \leq s K^2(s, 3, X) \; \; \text{ for } \; \; K \geq 2.
\end{equation}

\end{lemma}

A proof of Lemma \ref{thres} is given in Subsection
\ref{thresproof}. 

The range for the tuning parameter $\lambda$ is now specified as follows. 
For some constant $\frac{4 K(s, s, k_0)}{\eta} \leq M \leq
\frac{\sqrt{\Lambda_{\min}(s)}}{(1-\eta) c_0 \sigma}
\sqrt{\frac{n}{2\log p}}$, where  $0 < \eta < 1$, 
$\lambda_n$ is chosen such that
\begin{equation}
\label{eq::penalty-choice-instance}
16 M K(s, s, k_0) \geq 
\frac{\lambda_n}{c_0 \sigma \lambda_{\init} \sqrt{\bar{s}}}
\sqrt{\frac{n}{2 \log(p-s)}}
\geq \frac{64 K^2(s, s, k_0)}{\eta},
\end{equation}
where $\lambda_{\init}$ is defined in~\eqref{eq::init-lower-bound}
with $B = \sqrt{24}$ 
and $c_0 \geq 1$ is a small constant to be specified. 
The following theorem is an immediate result when we 
substitute $\tilde{\delta}_{\Sc}$ and 
$\tilde{\delta}_{S}$ that appear in Theorem~\ref{thm:first} 
with what we derived in Lemma~\ref{lemma:weights-bound}.

\begin{theorem}{\textnormal{(\bf{Variable selection for fixed design})}}
\label{thm:adaptive-recovery}
Consider the linear model in (\ref{eq::linear-model}) with fixed design
$X$, where $n \le p$, and each column of $X$ has its $\ell_2$-norm upper 
bounded by $\sqrt{n}$.
Suppose condition $RE(s, s, 3, X)$ (Assumption \ref{def:BRT-cond-two})
holds. Suppose on $\T$, for some $1> \eta > 0$,
$\lambda_n$ is chosen as in~\eqref{eq::penalty-choice-instance}
with $K(s, s, k_0) =  K(s, s, 3, X)$ and $c_0 =1$. 
Suppose $s$ satisfies~\eqref{eq::linear-sparsity} and 
\begin{eqnarray}
\label{eq::sparsity-condition}
\tilde{r_n} \sqrt{s} & \leq & \frac{1 - \eta}{32 K^2 \lambda_{\init}}, 
\; \; \text{ and }\\
\label{eq::tragic}
\beta_{\min}
& > &
\max\left\{ 
\frac{2 \tilde{r_n}}{1-\eta}, 
\frac{M}{\sqrt{3}}\right\} {16 K^2 \lambda_{\init} \sqrt{s}}
\end{eqnarray}
where $\lambda_{\init}$ is defined in~\eqref{eq::init-lower-bound}
with $B = \sqrt{24}$ and $K = K(s, s, 3, X)$.
Then, with probability $1 - 2/p^2$, the adaptive estimator 
in \eqref{eq::weighted-lasso} satisfies 
${\rm  supp}(\hat{\beta}) = {\rm supp}(\beta)$. 
\end{theorem}
A proof is given in Section \ref{sec:main-theorem-det}.
A first corollary follows immediately from 
Theorem~\ref{thm:adaptive-recovery} when we substitute
$\tilde{r_n} = \frac{\sqrt{s}}{\sqrt{\Lambda_{\min}(s)}}$ as shown
in Lemma \ref{lemma:norm-bounds}, formula ~\eqref{eq::r_n_bound-lemma} with
$c_0 =1$. 
\begin{corollary}{\textnormal{(\bf{Variable selection for fixed design:
 general bound for $\tilde r_n$})}}
\label{thm:adaptive-recovery-two}
Consider the linear model in (\ref{eq::linear-model}) with fixed design
$X$, where $n \le p$, and each column of $X$ has its $\ell_2$-norm upper 
bounded by $\sqrt{n}$.
Suppose that condition $RE(s, s, 3, X)$ (Assumption \ref{def:BRT-cond-two})
holds. Suppose that on $\T$ and for some $1> \eta > 0$, 
$\lambda_n$ is chosen as in~\eqref{eq::penalty-choice-instance}
with $K(s, s, k_0) =  K(s, s, 3, X)$ and $c_0 =1$,
\begin{eqnarray}
\label{eq::sparsity-condition-two}
s & \leq & 
\frac{\sqrt{\Lambda_{\min}(s)}(1-\eta)}{32 K^2 \lambda_{\init}}
\; \; \text{ and }  \; \; \; \\
\label{eq::tragic-instance}
\beta_{\min}
& > &
\max\left\{ 
\frac{2 \sqrt{s}}{(1-\eta)\sqrt{\Lambda_{\min}(s)}}, 
\frac{M}{\sqrt{3}}\right\} {16 K^2 \lambda_{\init} \sqrt{s}}
\end{eqnarray}
where $\lambda_{\init}$ is defined in~\eqref{eq::init-lower-bound}
with $B = \sqrt{24}$ and $K = K(s, s, 3, X)$.
Then, with probability $1 - 2/p^2$, the adaptive
estimator in \eqref{eq::weighted-lasso} satisfies 
${\rm supp}(\hat{\beta}) = {\rm supp}(\beta)$. 
\end{corollary}

Using the different bound $\tilde{r_n} = \frac{\theta_{s,
    1}\sqrt{s}}{\Lambda_{\min}(s)}$ from Lemma~\ref{lemma:norm-bounds}, formula
(\ref{eq::r_n_bound-lemma-2}), our next corollary shows that under
Assumption~\ref{def:BRT-cond-three},  
we can essentially achieve the sublinear sparsity level 
of~\eqref{eq::linear-sparsity} while conducting model selection.

\begin{corollary}{\textnormal{(\bf{Variable selection for fixed design:
special bound for $\tilde r_n$})}}
\label{thm:adaptive-recovery-three}
Consider the linear model in (\ref{eq::linear-model}) with fixed design
$X$, where $n \le p$, and each column of $X$ has 
$\ell_2$-norm upper bounded by $\sqrt{n}$. 
Suppose that Assumption~\ref{def:BRT-cond-three} holds
for $k_0 = 3$ and $m = s$.
Suppose that on $\T$ and for some $1> \eta > 0$, $\lambda_n$ is chosen as 
in~\eqref{eq::penalty-choice-instance}
with $K(s, s, k_0) =  K(s, s, 3, X)$ and $c_0 =1$. 
Suppose $s$ satisfies~\eqref{eq::linear-sparsity} and
\begin{eqnarray}
\label{eq::tragic-three}
\beta_{\min} & > &
\max\left\{ 
\frac{2 \sqrt{s} \theta_{1, s}}{(1-\eta)\Lambda_{\min}(s)}, 
 \frac{M}{\sqrt{3}} \right\}
{16 K^2 \lambda_{\init} \sqrt{s}}
\end{eqnarray}
where $\lambda_{\init}$ is defined in~\eqref{eq::init-lower-bound}
with $B = \sqrt{24}$ and $K = K(s, s, 3, X)$.
Then, with probability $1 - 2/p^2$, 
the adaptive estimator in \eqref{eq::weighted-lasso} satisfies 
${\rm supp}(\hat{\beta}) = {\rm supp}(\beta)$.
\end{corollary}
It is an open question whether the adaptive Lasso procedure can achieve
model selection consistency under such sparsity level under 
Assumption~\ref{def:BRT-cond-two} alone.

\section{The adaptive Lasso with random design}\label{adaprandom}
For a random design $X$ as in (\ref{eq::rand-des}), we make the following
assumption on $\Sigma$. 

\begin{assumption}
\textnormal{\bf{Restricted eigenvalue assumption} $RE(s, m, k_0, \Sigma)$}
\label{def:BRT-cond-random}
For some integer $1\leq s \leq p/2$, $m \geq s, s + m \leq p$, and
a positive number $k_0$, the following condition holds:
\beq
\label{eq::admissible-random}
\inv{K(s, m, k_0, \Sigma)} := \min_{\stackrel{J_0 \subseteq \{1, \ldots,
    p\},}{|J_0| \leq s}} 
\min_{\stackrel{\gamma \not=0,}
{\norm{\gamma_{J_0^c}}_1 \leq k_0 \norm{\gamma_{J_0}}_1}}
\; \;  \frac{\norm{\Sigma^{1/2} \gamma}_2}{\norm{\gamma_{J_{0m}}}_2} > 0.
\eeq
Suppose
\eqref{eq::eigen-cond-random} hold
and $\Sigma_{jj} = 1, \forall j = 1, \ldots, p$.
\end{assumption}
It is clear that in~\eqref{eq::admissible-random}, 
the set of admissible $\gamma$ is a superset of that 
in~\eqref{eq::eigen-admissible-random}. Hence we have:

\begin{proposition}
\label{prop:eigen-K-bound-random}
Suppose Assumption $RE(s,s, k_0, \Sigma)$ 
holds for some $1\leq s \leq p$ and some $k_0>0$. Then
$\frac{17\Lambda_{\min}(s)}{16} \geq \inv{ K^2(s, s, k_0, \Sigma)}$
for $\Lambda_{\min}(s)$ as defined in \eqref{eq::eigen-admissible-random}.
\end{proposition}

We now show that with high probability, 
Assumption $RE(s, m, k_0, X)$ holds for a random realization of 
$X$ whose row are i.i.d. vectors from $\sim N(0, \Sigma)$,
under Assumption~\ref{def:BRT-cond-random}, 
if $s = o\left(\sqrt{\frac{n}{\log p}}\right)$.

\begin{proposition}
\label{prop:bad-event}
Consider a random design $X$ as in (\ref{eq::rand-des}).
Assume that $\Sigma$ satisfies~\eqref{eq::admissible-random}. Then, on the
set $\X$ as defined in~\eqref{eq::good-random-design} and with $C_2$ as
in~\eqref{eq::good-random-design}, $X$ satisfies $RE(s, s, k_0, X)$ as in 
Assumption~\ref{def:BRT-cond-two}, with 
\begin{equation}
\label{eq::KX-upper-4}
K(s, s, k_0, X) \leq \sqrt{2} K(s, s, k_0, \Sigma),
\text{ for } 
s \leq \frac{\sqrt{{n}/{\log p}}}{32 C_2 K^2(s, 3, 3, \Sigma)}
\end{equation}
\end{proposition}
Its proof appears in Subsection~\ref{subsec:prop-bad-event}.

We can now state the result for a random design under 
Assumption~\ref{def:BRT-cond-random}.

\begin{theorem}{\textnormal{(\bf{Variable selection for a random design})}}
\label{thm:random-recovery}
Consider the linear model in (\ref{eq::linear-model}) with random design
$X$ as in (\ref{eq::rand-des}) with $n \le p$ and
 $p < e^{n/4C_2^2}$, where $C_2 > 4 \sqrt{5/3}$.
Suppose that Assumption \ref{def:BRT-cond-random} 
holds with $m = s$ and $k_0 = 3$. 
Suppose that on the set $\X \cap \T$ and for some $0< \eta < 1$, 
$\lambda_n$ is chosen as in~\eqref{eq::penalty-choice-instance} with 
$K(s, s, k_0) =  \sqrt{2} K(s, s, 3, \Sigma)$ and $c_0 =\sqrt{3/2}$;
suppose that
\begin{eqnarray}
\label{eq::sparsity-condition-random}
s \leq  \inv{32K^2(s, s, 3, \Sigma)}
\min\left\{\inv{C_2},
\frac{\sqrt{\Lambda_{\min}(s)} (1-\eta)}
{6 \sqrt{6} \sigma}
\right \}
\sqrt{\frac{n}{\log p}}
\end{eqnarray}
where $C_2$ is defined in~\eqref{eq::good-random-design}
In addition $\beta_{\min}$ satisfies~\eqref{eq::tragic-instance}
with $K = \sqrt{2} K(s, s, 3, \Sigma)$.
Then, with probability $1 - 3/p^2$, 
the adaptive Lasso estimator in \eqref{eq::weighted-lasso} satisfies 
${\rm supp}(\hat{\beta}) = {\rm supp}(\beta)$. 
\end{theorem}
A proof is given in Section \ref{subsec.thm-rand-recov}. 

\section{The adaptive Lasso in Gaussian graphical modeling}
\label{sec:graph}
Consider the problem of covariance selection described in Section
\ref{subsec.covsel}.
\subsection{The many regressions pursuit procedure}
The procedure for covariance selection in a Gaussian
graphical model based on a pursuit of many regressions has been proposed
and studied in \cite{MB06}. 

Consider $X_1,\ldots ,X_p \sim {\cal N}(0,\Sigma)$ as in
(\ref{eq:multGauss}). We can regress $X_i$ versus the other
variables $\{X_k;\ k \neq i\}$:
\begin{eqnarray}\label{eq::regr}
X_i = \sum_{j \neq i} \beta_j^i X_i + V_i
\end{eqnarray}
where $V_i$ is a normally distributed random variable with mean zero. Then,
denoting by $Q = \Sigma^{-1}$, it is well known that 
\begin{eqnarray}\label{eq::betaQ}
\beta_j^i = -\frac{Q_{ij}}{Q_{ii}}.
\end{eqnarray}
In particular, this implies that 
\begin{eqnarray*}
& &\mbox{there is an undirected edge between $i$ and $j$}\\
&\ \Leftrightarrow&\ \
\Sigma^{-1}_{ij} \neq 0\ 
\Leftrightarrow\ \beta_j^i \neq 0\ \mbox{and/or}\ \beta_i^j \neq 0,
\end{eqnarray*}
where the last statement holds due to the symmetry of $\Sigma^{-1}$.

The estimation of the edge set can then be done by one of the following
rules:
\begin{eqnarray*}
& &\mbox{there is an edge between $i$ and $j$}\ \Leftrightarrow\ \hat{\beta}_j^i \neq 0\
\mbox{and}\ 
\hat{\beta}_i^j \neq 0,\\
& &\mbox{there is an edge between $i$ and $j$}\ \Leftrightarrow\ \hat{\beta}_j^i \neq 0\
\mbox{or}\ 
\hat{\beta}_i^j \neq 0.
\end{eqnarray*}
Our obvious proposal is to use the adaptive Lasso estimates
$\hat{\beta}_{j;n}^i$ in the corresponding regressions as described in
(\ref{eq::regr}). The discrepancy between the ``and'' or
``or'' rule above vanishes with high probability.  

The theoretical analysis follows by our result for random
design linear models (Theorem \ref{thm:random-recovery}) and controlling
the error over $p$ different regressions. 
%
Let $\beta_{\min} = \min_{i, j} |\beta^{i}_j|$ and $s$ be the largest
node degree.
Our conditions on sparsity and $\beta_{\min}$ for linear models need to
hold for all $p$ regressions simultaneously and they are as follows.

\begin{assumption}\label{ass:gr1}
$\beta^i_j$ from (\ref{eq::regr}) satisfy the conditions on 
$\beta_{\min}$ as in~\eqref{eq::tragic-instance} 
$\forall i, j \in \{1, \ldots, p\}$ 
under Assumption \ref{def:BRT-cond-random}.
\end{assumption}
Equivalently, by
assuming $\Sigma_{jj}^{-1} = 1$ for all $j=1,\ldots ,p$ (see Assumption
\ref{def:BRT-cond-random}) and due to (\ref{eq::betaQ}), the non-zero
elements of $|\Sigma^{-1}_{ij}|$ are required to be upper-bounded by the
value of $\beta_{\min}$. 

\begin{assumption}
\label{ass:gr2}
The covariance matrix $\Sigma$ satisfies the restricted eigenvalue
condition in Assumption \ref{def:BRT-cond-random}.
In addition,~\eqref{eq::eigen-cond-random} is required to hold
on every subset $S \subset \{1, \ldots, p\}$ such that $|S| \leq s$.
\end{assumption}

\begin{assumption}\label{ass:gr4}
The size of the neighborhood, for all nodes, is bounded by an integer 
$1< s < p/2$ that satisfies~\eqref{eq::sparsity-condition-random}
under Assumption~\ref{def:BRT-cond-random}.
\end{assumption}

The following result can then be immediately derived using the union bound
for the $p$ regression in the many regressions pursuit. 

\begin{theorem}
\label{thm:mult-regres}
{\textnormal{(\bf{Covariance selection in Gaussian Graphical Models})}}
Consider the Gaussian graphical model with $n$
i.i.d. samples from (\ref{eq:multGauss}), where
$n \le p < e^{n/4C_2^2}$, where $C_2 > 4 \sqrt{5/3}$.
Suppose that Assumptions \ref{ass:gr1} -
\ref{ass:gr4} hold. Then,
\begin{eqnarray*}
\mathbb{P}\left({\rm supp}(\hat\Sigma^{-1}_n) = {\rm
    supp}(\Sigma^{-1})\right) \ge 1 - 3/p.
\end{eqnarray*}
\end{theorem}

\section{Discussion}\label{sec.discussion}
We have presented results for high-dimensional model selection in
regression and Gaussian graphical modeling. We make some
assumptions on (fixed or random) designs in terms of restricted
eigenvalues. Such assumptions are among the weakest for deriving oracle
inequalities in terms of $\|\hat{\beta} - \beta\|_q\ (q = 1,2)$
\citep{BRT08}. We show here that under such restricted eigenvalue
assumptions, the two-stage adaptive Lasso is able to correctly infer the
relevant variables in regression or the edge set in a Gaussian graphical
model. The ordinary Lasso can easily fail since the neighborhood stability
condition, or the equivalent irrepresentable condition, are necessary and
sufficient \citep{MB06,ZY07}. It is easy to construct examples
where the neighborhood stability condition fails but the restricted
eigenvalue condition holds for the situation where $n > p$, see for
example \cite{Zou06}. 

In the high-dimensional context, the relation between
the neighborhood stability condition and the restricted eigenvalue
assumption is not clear. However, the latter is a condition on an average
behavior (as an eigenvalue condition) while the former requires a relation
for a maximum: thus, we conjecture that the restricted eigenvalue assumption
is in general less restrictive than the neighborhood stability condition. 
In particular, although it appears non-trivial to derive a general
relation between these two conditions, one can certainly derive relations 
between them under additional assumptions; A thorough exposition of such 
relations is an interesting direction for future work, given the frequent
appearance of both types of conditions in the literature, 
for example in~\cite{MB06,ZY07,Wai08,CT07,MY08,BRT08}.
For high-dimensional Gaussian graphical modeling, using the reasoning
above, the restricted eigenvalue assumptions we make appears in general 
less restrictive (and easier to check) than the assumptions in \cite{MB06} 
and in \cite{RWRY08} who analyze the GLasso algorithm \cite{BGD07,FHT07}.

 
\section{Analysis of the weighted Lasso}\label{sec:sparsistence}
In the sequel, for clarity, we denote by $\beta^*$ the true parameter in
the linear model (\ref{eq::linear-model}). 
Inspired by the adaptive Lasso estimator defined 
in \eqref{eq::weighted-lasso}, we consider here the weighted Lasso with 
weights $0< w_j \ (j=1,\ldots ,p)$ which
solves the following optimization problem:  
\beq
\label{eq:weighted-lasso-spars}
\min_{\beta \in \R^p} 
\frac{1}{2n}\|Y-X\beta\|_2^2 + \lambda_n \sum_{j=1}^p w_j |\beta_j|,
\eeq
The only distinction between the adaptive and weighted Lasso is that we
assume that the weights are estimated in the former and pre-specified in
the latter approach. However, our theory below though does not 
depend whether the weights are random or not.  
For convenience we denote by 
\beq
\label{eq::max-min-def}
w_{\max}(S) = \max_{i \in S} \; w_{i}, \hspace{1cm}
w_{\min}(S^c) = \min_{j \in S^c} \; w_{j}.
\eeq

A slightly stronger notion than inferring the support of $\beta^*$ is the
recovery of the sign-pattern:
\begin{eqnarray*}
\sign(\hat{\beta}_n) = \sign(\beta^*).
\end{eqnarray*}
Furthermore, there are generally multiple solutions of the adaptive Lasso
estimator in (\ref{eq::weighted-lasso}) and in the weighted Lasso in
(\ref{eq:weighted-lasso-spars}). However, with high probability, the solution
is unique, see also Remark \ref{rem:unique} and Section
\ref{subsec.uniqueness}. 

As before, we denote by $\|A\|_{\infty} $ $= \max_{1 \le i \le k} \sum_{j=1}^m
|A_{ij}|$ for a $k \times m$ matrix $A$. 
First, let us state the following conditions that are imposed on the 
design matrix for the ordinary Lasso by~\cite{ZY07} and~\cite{Wai08}:
\begin{subeqnarray}
\label{eq:incoa}
\left\| X^T_{\Sc} X_S (X_S^T X_S)^{-1}\right\|_\infty & \leq & 1-\eta,
\;\;\text{for some $\eta\in(0,1]$, and} \\
\label{eq:incob}
\Lambda_{\min} \left(\onen X_S^T X_S\right) & \geq & \Lambda_{\min}(s) > 0,
\end{subeqnarray}
where $\Lambda_{\text{min}}(A)$ is the smallest eigenvalue of $A$.
Note that the second condition coincides with ours in~\eqref{eq::eigen-cond}.
\cite{MB06} formulated such conditions for a random design.

We impose the following incoherence conditions on the weighted Lasso.

\begin{definition}\textnormal{\bf{($(\vec{w}, S)$-incoherence condition)}}
\label{def:incoh-cond}
Let $X$ be an $n \times p$ matrix and let 
$S\subset \{1,\ldots, p\}$ be nonempty.
Let $\vec{w} = (w_1, w_2, \ldots, w_p)^T$ be a weight vector, where $w_j > 0
\forall j$. 
Let  $\vec{b}= (\sign(\beta^*_i) w_i)_{i \in S}$. 
We say that $X$ is \textit{$(\vec{w}, S)$-incoherent} if
for some $\eta\in(0,1),$
\begin{subeqnarray}
\label{eq:eta}
\forall j \in \Sc, \; \hspace{.5cm}
\abs{X^T_j X_S (X_S^T X_S)^{-1} \vec{b}} & \leq & w_j (1-\eta), \\
\label{eq:etb}
\Lambda_{\min}\left(\onen X_S^T X_S\right) & \geq & \Lambda_{\min}(s) > 0,
\end{subeqnarray}
where a sufficient condition for~\eqref{eq:eta} is
\begin{eqnarray}
\label{eq:eta-stronger}
\forall j \in \Sc, \hspace{1cm}
\norm{X^T_{\Sc} X_S (X_S^T X_S)^{-1}}_{\infty} & \leq & 
\frac{w_{\min}(\Sc)}{w_{\max}(S)} (1-\eta).
\end{eqnarray}
\end{definition}

We now state a general lemma about recovering the signs for the weighted Lasso
estimator as defined in~\eqref{eq::weighted-lasso}. 

\begin{lemma}{\textnormal{\bf(Sign recovery Lemma)}}
\label{lemma:recovery}
Consider the linear model in (\ref{eq::linear-model}) where the design
matrix $X$ satisfies~\eqref{eq:eta} and~\eqref{eq:etb}.
Let $c_0 = \max_{j \in \Sc} {\twonorm{X_j}}/{\sqrt{n}}$.
Suppose that $w_j > 0, \forall j=1,\ldots ,p$ and $\lambda_n$ is chosen
such that 
\ben
\label{eq:thm-cond-lambda}
\lambda_n w_{\min}(\Sc) 
& \geq &
 \frac{4 c_0 \sigma}{\eta}\sqrt{\frac{2\log(p-s)}{n}},
\een
where $w_{\min}(\Sc), w_{\max}(S)$ are as defined
in~\eqref{eq::max-min-def}. Furthermore, assume
\begin{eqnarray}
\label{eq::beta-min-recovery}
\beta_{\min} 
& > &
\max\left\{
\frac{4 c_0 \sigma}{\Lambda_{\min}(s)} \sqrt{\frac{6 s \log p}{n}},
\frac{2 \lambda_n w_{\max}(S) \sqrt{s}}{\Lambda_{\min}(s)} \right\}
\end{eqnarray}
Then for $\hat \beta$ in (\ref{eq:weighted-lasso-spars}):
\begin{equation*} 
\label{eq::perfect-recovery}
\prob{\sign(\hat\beta) = \sign(\beta^*)}
\geq 1- 2/p^2,
\end{equation*}
Moreover, 
with $\T$ defined in~\eqref{eq::low-noise},
we have
$\prob{(\sign(\hat\beta) \not= \sign(\beta^*)) \cap \T} 
\leq 2/p^2$.
\end{lemma}
A proof is given in Section~\ref{subsec.lemma-recovery}. 
Note that in case $w_{\min}(\Sc) = w_{\max}(S) = 1$,
conditions~\eqref{eq:eta-stronger} and~\eqref{eq::beta-min-recovery} 
reduce to~\eqref{eq:incoa} and the the statement of Lemma
\ref{lemma:recovery} is exactly the same 
as Theorem~$1$ in~\cite{Wai08}. 

\section{Proof of Lemma~\ref{lemma:twonorm-random}}
\label{subsec:lemmatwonorm}    
  
\begin{lemma}
\label{lemma:gaussian-noise}
For fixed design $X$ with $\max_j \|X_j\|_2 \le c_o \sqrt{n}$ we have
for $\T$ as defined in~\eqref{eq::low-noise},
\beq
\prob{\T^c} \leq 1/p^2.
\eeq
\end{lemma}

\begin{proof}
Define the random variables
$$Y_j = \inv{n} \sum_{i=1}^n \e_i X_{i,j}.$$
Note that $\max_{1 \le j \le p} |Y_j| = \|X^T \epsilon/n\|_{\infty}$. 
We have $\expct{Y_j} = 0$ and
$\var(Y_j) = \frac{\twonorm{X_j}^2\sigma_{\e}^2}{n^2} 
\leq \frac{c_0 \sigma_{\e}^2}{n}$. 
Obviously, $Y_j$ has its tail probability dominated by that of
$Z \sim N(0, \frac{c_0^2 \sigma_{\e}^2}{n})$:
\begin{eqnarray*}
\prob{|Y_j| \geq t} \leq  \prob{|Z| \geq t} \leq 
\frac{c_0 \sigma_{\e}}{\sqrt{n} t} \exp\left(\frac{-n t^2}{2 c_0^2
    \sigma_{\e}^2}\right).
\end{eqnarray*}
We can now apply the union bound to obtain:
\ben
\prob{\max_{1 \leq j \leq p} |Y_j| \geq t } &\leq& 
p  \frac{c_0 \sigma_{\e}}{\sqrt{n} t} \exp\left(\frac{-n t^2}{2 c_0^2
    \sigma_{\e}^2}\right) \\
&=& 
\exp\left(-\left(\frac{n t^2}{2 c_0^2 \sigma_{\e}^2} + 
\log \frac{t \sqrt{n}}{c_0 \sigma_{\e}} - \log p\right) \right).
\een
By choosing $t = c_o \sigma_{\e} \sqrt{6 \log(p)/n}$, the right-hand side
is bounded by $1/p^2$.
\end{proof}


We now show that $\prob{\X} \ge 1 - 1/p^2$. 

We denote $\Sigma_{ii} := \sigma_i^2$ 
throughout the rest of this proof.
We first state the following large inequality bound for the nondiagonal 
entries of $\Sigma$, adapted from Lemma 38~\citep{ZLW08} by plugging in 
$\sigma_{i}^2 = 1, \forall i=1, \ldots, p$ and using the 
fact that 
$|\Sigma_{jk}| = |\rho_{jk} \sigma_{j} \sigma_k| \leq 1, \forall j \not= k$,
where $\rho_{jk}$ is the correlation coefficient between variables
$X_j$ and $X_k$.

\begin{lemma}{\textnormal{\citep{ZLW08}}}
\label{lemma:boxcar-deviation-non-diag}
Let $\Psi_{jk} = (1 + \Sigma^2_{jk})/2$. For $0 \leq \tau \leq \Psi_{jk}$, 
\begin{eqnarray}
\label{eq::boxcar-non-diag}
\prob{|\Delta_{jk}| > \tau} \leq 
\exp\left\{-\frac{3 n \tau^2}{10(1 + \Sigma_{jk}^2)} \right\} 
\leq
\exp\left\{-\frac{3 n \tau^2}{20} \right\}.
\end{eqnarray}
\end{lemma}
We now also state a large deviation bound for the $\chi^2_n$ 
distribution~\cite{Joh01}:
\begin{eqnarray}
\label{eq::chi-dev}
\prob{\frac{\chi^2_n}{n} - 1 > \tau} & \leq & 
\exp\left(\frac{-3 n \tau^2}{16}\right),\; \text{for} \;0 \leq \tau \leq \half.
\end{eqnarray}
Hence by the union bound, we have $j = 1, \ldots, p$, for $\tau < 1/2$,
\begin{eqnarray}
\label{eq::chi-dev-union}
\prob{\max_{j=1, \ldots, p} \frac{\|X_j\|_2^2}{n}   - 1 >  \tau}
\leq p \exp\left(\frac{-3 n \tau^2}{16}\right). 
\end{eqnarray}

\begin{lemma}
\label{lemma:gaussian-covars} 
For a random design $X$ as in (\ref{eq::rand-des}) with 
$\Sigma_{jj} = 1, \forall j \in \{1, \ldots, p\}$, 
and for $p < e^{n/4 C_2^2}$, where $C_2 > 4 \sqrt{5/3}$,
we have 
\begin{eqnarray*}
\prob{\X} \ge 1 - 1/p^2.
\end{eqnarray*}
\end{lemma}

\begin{proof}
Now it is clear that we have $p(p-1)/2$ unique non-diagonal
entries $\sigma_{jk}, \forall j \not= k$ and $p$ diagonal entries. 
By the union bound and by taking
$\tau =  C_2 \sqrt{\frac{\log p}{n}}$ in~\eqref{eq::chi-dev-union} 
and~\eqref{eq::boxcar-non-diag}, we have
\label{eq::chi-dev}
\begin{eqnarray*}
\label{eq::Delta-bound}
\prob{\X^c} 
& = & 
\prob{\max_{jk} |\Delta_{jk}| \geq C_2 \sqrt{\frac{\log p}{n}}} \\
& \leq & 
p \exp\left(- \frac{3 C_2^2 \log p}{16}\right) + 
\frac{p^2 - p}{2} \exp\left(- \frac{3 C_2^2 \log p}{20}\right) \\
& \leq &
p^2 \exp\left(- \frac{3 C_2^2 \log p}{20}\right)
= p^{- \frac{3C_2^2}{20} + 2}  <  \frac{1}{p^{2}}
\end{eqnarray*}
for $C_2 > 4 \sqrt{5/3}$.
Finally, $p < e^{n/4 C_2^2}$ guarantees that 
$C_2 \sqrt{\frac{\log p}{n}} < 1/2$.
\end{proof}


\section{Proofs for Section \ref{adapfix}}

Throughout this section, we have 
$\lambda_{\init} = B c_0 \sigma_{\e} \sqrt{\frac{\log p}{n}}$ 
with $B = \sqrt{24}$.

\subsection{The Lasso as initial estimator}
\label{sec:main-result}
Lemma~\ref{lemma:weights-bound} crucially uses 
the bound on the $\ell_1$-loss of the initial Lasso estimator.

Our proof follows that of~\cite{BRT08}.
Let $\beta_{\init}$ be as in~\eqref{eq::initial-estimator} and $\delta =
\beta_{\init} - \beta^*$. 
The set $\T$ is defined in~\eqref{eq::low-noise}.
We first show Lemma~\ref{lemma:magic-number}; we then apply condition 
$RE(s, k_0, X)$ on $\delta$ with $k_0 = 3$ under $\T$ to
derive various norm bounds.

\begin{lemma}
\label{lemma:magic-number}
For fixed design, on ${\T}$, $\norm{\delta_{S^c}}_1 \leq
3\norm{\delta_{S}}_1$. 
\end{lemma}
\begin{proof}
Since $\beta_{\init}$ is a Lasso solution, we have
\begin{eqnarray*}
\lambda_{\init} \norm{\beta^*}_1 - 
\lambda_{\init} \norm{\beta_{\init}}_1 
& \geq & 
\inv{2n} \norm{Y- X  \beta_{\init}}^2_2  - 
\inv{2n} \norm{Y- X \beta^* }^2_2 \\
& \geq & 
\inv{2n} \| X \delta \|^2_2 - \frac{\delta^T X^T \e}{n} 
\end{eqnarray*}
Hence on the set $\T$ as in~\eqref{eq::low-noise}, we have
\begin{eqnarray}
 \nonumber
\norm{X \delta}_n^2 
& \leq &
2 \lambda_{\init} \norm{\beta^*}_1
-  2 \lambda_{\init} \norm{\beta_{\init}}_1 +  
2 \norm{\frac{X^T \e}{n}}_{\infty} \norm{\delta}_1 \\
\label{eq::precondition}
& \leq & \lambda_{\init} \left(2\norm{\beta^*}_1 
- 2 \norm{\beta_{\init}}_1 + \norm{\delta}_1\right),
\end{eqnarray}
where by the triangle inequality, and $\beta^*_{\Sc} = 0$, we have
\begin{eqnarray} 
\nonumber
0 & \leq & 2 \norm{\beta^*}_1 -  2 \norm{\beta_{\init}}_1 + \norm{\delta}_1 
\\ \nonumber
& = &
2 \norm{\beta^*_S}_1 - 2 \norm{\beta_{S, \init}}_1 -
2 \norm{\delta_{\Sc}}_1 + \norm{\delta_S}_1 + \norm{\delta_{S^c}}_1  \\
\label{eq::magic-number-2}
& \leq &
3 \norm{\delta_{S}}_1 - \norm{\delta_{S^c}}_1.
\end{eqnarray}
Thus Lemma~\ref{lemma:magic-number} holds.
\end{proof}

\begin{proposition}
{\textnormal{({\bf $\ell_p$-loss for the initial estimator}, \citep{BRT08})}}
\label{prop:initial-bound}
Consider the linear model in (\ref{eq::linear-model}) with fixed design
satisfying $\max_j \twonorm{X_j} \leq c_0 \sqrt{n}$. Suppose that $RE(s, 3, X)$
holds. Let $\delta = \beta_{\init} - \beta^*$ with $\beta_{\init}$ defined
in~\eqref{eq::initial-estimator} with 
 $$\lambda_{\init} = B c_0 \sigma_{\e} \sqrt{\frac{\log p}{n}}.$$
Then, on the set $\T$ in~\eqref{eq::low-noise},
\begin{eqnarray}
\label{eq::magic-S}
\norm{\delta_S}_2 & \leq & 4 K^2(s, 3, X) \lambda_{\init} \sqrt{s}. \\
\label{eq::magic-thm}
\norm{\delta}_1 & \leq & 4 K^2(s, 3, X) \lambda_{\init} s;
\end{eqnarray}
Moreover, under the stronger assumption $RE(s, s, 3, X)$, and
on the set $\T$ as in~\eqref{eq::low-noise},
\begin{eqnarray}
\label{eq::BRT-2-bound}
\norm{\delta}_2 & \leq & 16 K^2(s, s, 3, X) \lambda_{\init} \sqrt{s}.
\end{eqnarray}
\end{proposition}

\begin{proof}
On the set $\T$, by~\eqref{eq::precondition}
and~\eqref{eq::magic-number-2},
\begin{eqnarray}
\nonumber
\norm{X \delta}_n^2 +  \lambda_{\init} \norm{\delta}_1
& \leq &
\lambda_{\min} 
\left(3 \norm{\delta_{S}}_1 - \norm{\delta_{S^c}}_1 + \norm{\delta_S}_1 
+ \norm{\delta_{\Sc}}_1 \right)  \\
\label{eq::last} 
& = &  
4 \lambda_{\init} \norm{\delta_S}_1
\leq 4\lambda_{\init} \sqrt {s}  \norm{\delta_S}_2 \\
\label{eq::last-2} 
& \leq &
4 \lambda_{\init} \sqrt{s} K(s, 3, X) \norm{X \delta}_n  \\ \nonumber
& \leq &
4 K^2(s, 3, X) \lambda_{\init}^2 s + \norm{X \delta}_n^2,
\end{eqnarray}
where~\eqref{eq::last-2} is due to condition $RE(s, 3, X)$ 
and Lemma~\ref{lemma:magic-number}. Hence~\eqref{eq::magic-thm} holds. 
Now by $RE(s, 3, X)$ and~\eqref{eq::last}, we have
\begin{eqnarray}
\norm{\delta_S}^2_2 \leq 
K^2(s, 3, X)\norm{X \delta}^2_n 
& \leq & 
K^2(s, 3, X) 4 \lambda_{\init} \sqrt {s}  \norm{\delta_S}_2.
\end{eqnarray}
Hence~\eqref{eq::magic-S} holds.
Finally, on the set $\T$, given Lemma~\ref{lemma:magic-number}, 
by $RE(s, s, 3, X)$ and~\eqref{eq::last},
we have
\ben
\norm{\delta_{S S'}}_2^2 & \leq &
K^2(s, s, 3, X) \norm{X \delta}_n^2 \\
& \leq & 
K^2(s, s, 3, X) 4 \lambda_{\init} \sqrt{s} \norm{\delta_S}_2 \\
& \leq &
K^2(s, s, 3, X) 4 \lambda_{\init} \sqrt{s} \norm{\delta_{SS'}}_2.
\een
Hence from the following inequality~\eqref{eq::BRT-SS-bound} (e.g.,
cf. (B.28) in~\cite{BRT08})
\begin{eqnarray}
\label{eq::BRT-SS-bound}
\norm{\delta}_2 & \leq & (1 + k_0)  \norm{\delta_{S S'}}_2,
\end{eqnarray}
we obtain~\eqref{eq::BRT-2-bound}.
\end{proof}

\subsection{
Proof of Lemma~\ref{lemma:weights-bound}}\label{sec:initial}
By Proposition~\ref{prop:initial-bound}, and (B.26) in~\cite{BRT08}, 
\ben 
\norm{\delta_S}_2 & \leq & 4 K(s, 3, X)^2 \lambda_{\init} \sqrt{s}, \\
\label{eq::total-sum} 
\norm{\delta}_1 & \leq & 4 K(s, 3, X)^2 \lambda_{\init} s, \; \; \text { where} \\
\norm{\delta_{\Sc}}_1 & \leq & 3 \norm{\delta}_1,
\een
due to a property of the Lasso estimator (see, for example~\cite{BRT08}).
This allows us to conclude that on the set $\T$ as in~\eqref{eq::low-noise}, 
\begin{eqnarray}
\label{eq::partial-sum} 
\norm{\delta_S}_{\infty} & \leq &
\norm{\delta_S}_{2} \leq 4 K(s, 3, X)^2 \lambda_{\init} \sqrt{s}, \\
\label{eq::partial-sum-2} 
\norm{\delta_{\Sc}}_1 & \leq & \frac{3}{4} \norm{\delta}_1 \leq 
3 K(s, 3, X)^2 \lambda_{\init} s.
\end{eqnarray}
Thus we have by~\eqref{eq::weight-bounds-conditions},
~\eqref{eq::partial-sum} and~\eqref{eq::partial-sum-2},
\begin{eqnarray}
\label{eq::S-init-wights}
\forall \ i \in S, \; \;
|\beta_{i, \init}| & \geq & \beta_{\min} - \norm{\delta_S}_{\infty} 
\geq 4 K(s, 3, X)^2 \lambda_{\init} \sqrt{s}, \\
\forall j \in \Sc, \; \;
|\beta_{j, {\init}} | & \leq &  \norm{\delta_{\Sc}}_{\infty} \leq
 \norm{\delta_{\Sc}}_1 \leq  3 K(s, 3, X)^2 \lambda_{\init} s.
\end{eqnarray}
\qed

\subsection{Proof of Lemma \ref{thres}}\label{thresproof}
If we threshold $\beta_{\init}$ at the value of 
$4 \lambda_{\init}$, by (\ref{eq::S-init-wights}), we have $\bar{S}
\supseteq S$. Moreover, by~\eqref{eq::partial-sum-2},
we include at most ${3K(s, 3, X)^2 s}/{4}$ more entries from 
$\Sc$ in $\bar{S}$; thus for $K(s, 3, X) \geq 2$,
$$s \leq |\bar{S}| \leq s + \frac{3 s K(s, 3, X)^2}{4} \leq s K(s, 3, X)^2.$$  
In addition, we have $\forall j \in \Sc$, by~\eqref{eq::BRT-2-bound},
\ben
|\beta_{j, {\init}} | & \leq & \norm{\delta_{\Sc}}_{\infty}
\leq \norm{\delta_{\Sc}}_{2} \\
& \leq & 16 K^2(s, s, 3, X) \lambda_{\init} \sqrt{s},
\een
under Assumption $RE(s, s, 3, X)$ and condition $\T$.
\qed

\subsection{Proof of Theorem~\ref{thm:adaptive-recovery}}
\label{sec:main-theorem-det}
It is clear that once we finish checking conditions on 
$\lambda_n$ in~\eqref{eq::penalty-choice}, on $s$ as 
in~\eqref{eq::sparsity-condition-first} and 
on $\beta_{\min}$ as in~\eqref{eq::tragic-0} hold,
we can invoke Theorem~\ref{thm:first} to finish the proof.
Formula~\eqref{eq::weight-bounds-conditions} is satisfied
assuming~\eqref{eq::tragic-instance}.  
Hence by choosing
\begin{eqnarray}
\label{eq::norm-S}
\tilde{\delta}_{S} & :=  &
4 K^2(s, 3, X) \lambda_{\init} \sqrt{s},\\
\label{eq::norm-Sc}
\tilde{\delta}_{\Sc} 
& := & 16 K^2(s, s, 3, X) \lambda_{\init} \sqrt{s},
\end{eqnarray}
we have
$\tilde{\delta}_{S} \geq \norm{\delta_{S}}_{\infty}$
and
$\tilde{\delta}_{\Sc} \geq \norm{\delta_{\Sc}}_{\infty}$
by~\eqref{eq::e1} and~\eqref{eq::weight-bounds-min-two}.
Now by~\eqref{eq::superset},
\begin{eqnarray*}
\lambda_n & \geq & 
\frac{64 \sigma K^2(s, s, 3, X) \lambda_{\init} \sqrt{\size{\bar{S}}}}{\eta} 
\sqrt{\frac{2 \log(p-s)}{n}} \\ \nonumber
& \geq & 
\frac{4 \sigma}{\eta}\sqrt{\frac{2\log(p-s)}{n}}
16 K^2(s, s, 3, X) \lambda_{\init} \sqrt{s} \\ \nonumber
& = & 
\frac{4 \sigma \tilde{\delta}_{\Sc}}{\eta}
\sqrt{\frac{2\log(p-s)}{n}}
\end{eqnarray*}
and
\begin{eqnarray*}
\lambda_n
& \leq &
\frac{16 M K^2(s, 3, 3, X) \sigma \lambda_{\init} \sqrt{\size{\bar{S}}}}
{K(s, 3, X)}\sqrt{\frac{2 \log(p-s)}{n}} \\
& \leq & 
M \sigma 16 K^2(s, 3, 3, X) \lambda_{\init}\sqrt{s}
\sqrt{\frac{2 \log(p-s)}{n}} \\
& = &
M \sigma \tilde{\delta}_{\Sc}
\sqrt{\frac{2 \log(p-s)}{n}},
\end{eqnarray*}
and thus~\eqref{eq::penalty-choice} holds with $c_0 = 1$.
Furthermore, for the sparsity $s$, \eqref{eq::sparsity-condition}
guarantees that \eqref{eq::sparsity-condition-first} holds
by~\eqref{eq::norm-Sc}. Finally, regarding $\beta_{\min}$,
\eqref{eq::tragic-0} holds given \eqref{eq::tragic}, 
as $\frac{16 M K^2 \lambda_{\init} \sqrt{s}}{\sqrt{3}}$ clearly dominates
the first and the third term in \eqref{eq::tragic-0} by the definition 
of~\eqref{eq::norm-S} and~\eqref{eq::norm-Sc}, and
the fact that $\inv{\Lambda_{\min}(s)} \leq K^2(s, k_0, X)$ by 
Proposition~\ref{prop:eigen-K-bound}; and it also dominates the second term 
given~\eqref{eq::linear-sparsity} and the upper bound on $\lambda_n$. 
\qed

\subsection{Bounds for $r_n$}

\begin{lemma}
\label{lemma:norm-bounds}
Consider a fixed design $X$ with $\max_j \twonorm{X_j} \leq c_0 \sqrt{n}$
and assume that~\eqref{eq::eigen-cond} holds. Then for all subsets $S$ with
$|S| \le s$,  
\begin{eqnarray}
\label{eq::r_n_bound-lemma}
r_n := \norm{X^T_{S^c} X_S (X_S^T X_S)^{-1}}_{\infty} 
 & \leq & \frac{c_0 \sqrt{s}}{\sqrt{\Lambda_{\min}(s)}}. \\
\label{eq::r_n_bound-lemma-2}
r_n & \leq & \frac{\theta_{1, s}\sqrt{s}}{\Lambda_{\min}(s)},
\end{eqnarray}
where $\theta_{1,s}$ is given in \ref{label:correlation-coefficient}
\end{lemma}

\begin{proof}
As a shorthand, we let $P_S = X_S (X_S^T X_S)^{-1} X_S^T$ denote the 
projection matrix and define
$$\forall j \in S^c, \; \; r_j =(X_S^T X_S)^{-1}X_S^T X_j.$$
Bounding $\norm{r_j}_1\ \forall j$ yields a bound on $r_n$.
First we have for all $j \in \Sc$,
\begin{eqnarray}
\label{eq::first-bound}
\twonorm{X_S r_j} = \twonorm{X_S (X_S^T X_S)^{-1} X_S^T X_j} = 
\twonorm{P_S X_j} \end{eqnarray}
$$
\leq \twonorm{X_j} \leq c_0 \sqrt{n}.
$$
On the other hand, by the restricted eigenvalue assumption, we have
$$\twonorm{X_S r_j}^2 = r_j^T X_S^T X_S r \geq
n \Lambda_{\min} \left(\frac{X_S^T X_S}{n} \right) \twonorm{r_j}^2.$$
Thus we have that $\twonorm{r_j} 
\leq \frac{c_0}{\sqrt{\Lambda_{\min}(s)}}, \forall j \in S^c$, and hence
$$r_n = \max_{j \in S^c} \norm{r_j}_1 \leq 
\max_{j \in S^c} \sqrt{s} \norm{r_j}_2 = \sqrt{s} \max_{j \in S^c} \norm{r_j}_2
\leq \frac{c_0 \sqrt{s}}{\sqrt{\Lambda_{\min}(s)}}.$$

Next we note that using~\eqref{label:correlation-coefficient}, we can
bound $r_n$ as follows, which has
essentially been shown in~\cite{CT07}.
For $P_S X_j = X_S r_j$, with
\begin{eqnarray*}
\twonorm{r_j} 
& \leq & \frac{\twonorm{X_S r_j}}{\sqrt{n \Lambda_{\min}(s)}} 
= \frac{\twonorm{P_S X_j}}{\sqrt{n \Lambda_{\min}(s)}}
\end{eqnarray*}
we have
\begin{eqnarray*}
\label{eq::second-bound}
\frac{\twonorm{P_S X_j}^2}{n} 
& = & \frac{\ip{P_S X_j, X_j}}{n}  = \frac{\ip{X_S r_j, X_j}}{n} \\ 
& \leq & \theta_{1,s} \twonorm{r_j} \leq 
\theta_{1,s} \frac{\twonorm{X_S r_j}}{\sqrt{n \Lambda_{\min(s)}}}  = 
\theta_{1,s} \frac{\twonorm{P_S X_j}}{\sqrt{n \Lambda_{\min}(s)}}
\end{eqnarray*}
Hence,
\begin{eqnarray*}
\twonorm{P_S X_j} & \leq & \frac{\sqrt{n}
  \theta_{1,s}}{\sqrt{\Lambda_{\min}(s)}}\ \ \mbox{and}\ \ r_n \leq \frac{
  \sqrt{s} \theta_{1,s}} {\Lambda_{\min}(s)}.
\end{eqnarray*}
\end{proof}

\section{Proofs for Section \ref{adaprandom}}

\subsection{Proof of
  Proposition~\ref{prop:bad-event}} \label{subsec:prop-bad-event} 
We first bound $\norm{X \gamma}_n^2 - \gamma^T \Sigma \gamma$.
\begin{eqnarray*}
& &\abs{\norm{X \gamma}_n^2 - \gamma^T \Sigma \gamma}  
= 
\abs{\gamma^T \hat\Sigma \gamma - \gamma^T \Sigma \gamma}
 =  
\abs{\sum_{j=1}^p \sum_{k=1}^p  \gamma_j \gamma_k (\hat\Sigma_{jk} - \Sigma_{jk})} \\
& \leq & 
\abs{\sum_{j \in S} \sum_{k \in S} \gamma_j \gamma_k (\hat\Sigma_{jk} - \Sigma_{jk})} 
+ 
\abs{\sum_{j \in S^c} \sum_{k \in S^c} \gamma_j \gamma_k (\hat\Sigma_{jk} -
  \Sigma_{jk})}\\
&+& 2 \abs{\sum_{j \in S} \sum_{k \in S^c} \gamma_j \gamma_k
  (\hat\Sigma_{jk} - \Sigma_{jk})} \\ 
& \leq & 
\max_{j, k} |\Delta_{jk}| \left(
\norm{\gamma_S}_1^2 + 2 \norm{\gamma_{S}}_1 \norm{\gamma_{S^c}}_1 +
\norm{\gamma_{S^c}}^2_1
\right),
\end{eqnarray*}
where $\Delta = \hat{\Sigma} - \Sigma$. 
Now given that $\norm{\gamma_{S^c}}_1 \leq k_0 \norm{\gamma_{S}}_1$
and $\norm{\gamma_S}^2_1 \leq s \twonorm{\gamma_S}^2$, we have 
\begin{eqnarray*}
& &\abs{\norm{X \gamma}_n^2 - \gamma^T \Sigma \gamma} 
\leq 
\max_{j, k} |\Delta_{jk}|  \norm{\gamma_S}^2_1 (1 + 2 k_0 + k_0^2)  \\
& \leq &
\max_{j, k} |\Delta_{jk}|  \norm{\gamma_S}^2_1 (1 + k_0)^2  \leq 
s (1 + k_0)^2
\max_{j, k} |\Delta_{jk}| \twonorm{\gamma_S}^2.
\end{eqnarray*}
Let $\gamma_{SS'} = \gamma_S \cup \gamma_{S'}$, where 
$\gamma_{S'}$ denote the subset of $\{1, \ldots, p\}$
corresponding to the $s$ largest  coordinates of $\gamma$ in their 
absolute values in $\gamma_{\Sc}$. We have on $\X$, using
Assumption~\ref{def:BRT-cond-random}, 
\begin{eqnarray*}
& &\norm{X \gamma}_n^2 
\geq  
\gamma^T \Sigma \gamma -
s (1 + k_0)^2 \max_{j, k} |\Delta_{jk}| \twonorm{\gamma_S}^2  \\
& \geq & 
\frac{\twonorm{\gamma_{SS'}}^2}{K(s, s, k_0, \Sigma)^2} -
s (1 + k_0)^2 \max_{j, k} |\Delta_{jk}| \twonorm{\gamma_S}^2 \geq  
\frac{\twonorm{\gamma_{SS'}}^2 }{2 K(s, s, k_0, \Sigma)^2},
\end{eqnarray*}
and hence~\eqref{eq::KX-upper-4} holds.
\qed

\subsection{Eigenvalue bounds}
We now show that~\eqref{eq::eigen-cond} is satisfied with high 
probability for a random design $X$, given its population correspondent as in 
~\eqref{eq::eigen-cond-random}.

\begin{lemma}
\label{lemma:norm-bounds-random}
Let $X$ be a random design as in (\ref{eq::rand-des}).
Let $s \leq \frac{\Lambda_{\min}(s)}{16 C_2} \sqrt{\frac{n}{\log p}}$ for
$C_2$ as defined in~\eqref{eq::good-random-design}.
We have on the set $\X$,
\begin{eqnarray}
\label{eq::eigen-exist}
\Lambda_{\min} \left(\frac{X_S^T X_S}{n}\right) \geq  \Lambda_{\min}(s),
\end{eqnarray}
for all subsets $S \subset \{1, \ldots, p\}$ with $|S| \leq s$ 
where~\eqref{eq::eigen-cond-random} hold.
\end{lemma}
 
\begin{proof}
On the set $\X$, for all subsets $S$ with $|S| \leq s$,
\begin{eqnarray}
\label{eq::matrix-perturb}
\abs{
\Lambda_{\min} \left(\frac{X_S^T X_S}{n}\right)
-\Lambda_{\min}(\Sigma_{SS})} 
& \leq & 
\twonorm{\left(\frac{X_S^T X_S}{n}\right) - \Sigma_{SS}} \\
\label{eq::symmetric}
& \leq &
 \norm{\left(\frac{X_S^T X_S}{n}\right) - \Sigma_{SS}}_{\infty} \\
\label{eq::cite-Delta-bound}
& \leq &  s C_2 \sqrt{\frac{\log p}{n}} 
\leq \frac{\Lambda_{\min}(s)}{16},
\end{eqnarray}
where $\|.\|_2$ denotes here the operator norm of a matrix. 
\eqref{eq::matrix-perturb} is a standard result in matrix perturbation
theory,~\eqref{eq::symmetric} is due to the fact that $\hat\Sigma$ and $\Sigma$
are symmetric, and~\eqref{eq::cite-Delta-bound} is due 
to~\eqref{eq::good-random-design} and the bound on $s$.
Hence for all subsets $S$ with $|S| \leq s$
that satisfy
$\Lambda_{\min} (\Sigma_{SS} ) \geq \frac{17}{16} \Lambda_{\min}(s)$ (as
defined in (\ref{eq::eigen-admissible-random})),
\eqref{eq::eigen-exist} holds.
\end{proof}

\subsection{Proof of Theorem~\ref{thm:random-recovery}} 
\label{subsec.thm-rand-recov} 
     
As corollary of Lemmas \ref{lemma:norm-bounds} and 
\ref{lemma:norm-bounds-random}, we  have

\begin{corollary}\label{iets}
\label{coro:r_n-bounds} Consider a random design $X$. Then on the set
$\X$ defined in~\eqref{eq::good-random-design},
~\eqref{eq::r_n_bound-lemma} holds with $c_0 =\sqrt{3/2}$, for all subsets $S$ with 
$|S| \le s$.
\end{corollary}

It is clear that~\eqref{eq::weight-bounds-conditions} is always
satisfied given~\eqref{eq::tragic-instance}, 
where $K = \sqrt{2} K(s, 3, 3, \Sigma)$,
as $K(s, s, k_0, X) \leq \sqrt{2} K(s, s, k_0, \Sigma)$
 by Proposition~\ref{prop:bad-event}.
We now show that the conditions on $\lambda_n$, $s$ and $\beta_{\min}$ as 
required by Theorem~\ref{thm:first} are satisfied on $\X \cap \T$.
First we take
\begin{eqnarray}
\label{eq::norm-S-random}
\tilde{\delta}_{S} & :=  &
8 K^2(s, s, 3, \Sigma) \lambda_{\init} \sqrt{s}  \\
\label{eq::norm-Sc-random}
\tilde{\delta}_{\Sc} 
& := &   32 K^2(s, s, 3, \Sigma) \lambda_{\init} \sqrt{s} , \\ 
\label{eq::r-n-upper}
\tilde{r_n} & := & \frac{\sqrt{3 s}}{\sqrt{2 \Lambda_{\min}(s)}},
\end{eqnarray}
where~\eqref{eq::r-n-upper} holds by Corollary~\ref{coro:r_n-bounds},
for which
$$s \leq \inv{32 C_2 K^2(s, s, 3, \Sigma)}
\leq \frac{\Lambda_{\min}(s)}{16 C_2} \sqrt{\frac{n}{\log p}},$$
by Proposition~\ref{prop:eigen-K-bound-random}. It is clear that 
\begin{eqnarray}
\tilde{\delta}_{S} & \geq &
4 K^2(s, s, 3, X) \lambda_{\init} \sqrt{s} \geq \norm{\delta_{S}}_{\infty},
 \text{ and } \\ 
\tilde{\delta}_{\Sc} 
& \geq & 16 K^2(s, s, 3, X) \lambda_{\init} \sqrt{s}
 \geq \norm{\delta_{\Sc}}_{\infty}
\end{eqnarray}
given~\eqref{eq::e1} and~\eqref{eq::weight-bounds-min}, and 
Proposition~\ref{prop:bad-event}.
Regarding the condition on $\lambda_n$,
by Proposition~\ref{prop:bad-event},~\eqref{eq::superset} 
and~\eqref{eq::norm-Sc-random}, we have
\begin{eqnarray*}
\lambda_n & \geq & 
\frac{128 c_0 \sigma K^2(s, s, 3, \Sigma) \lambda_{\init} \sqrt{\size{\bar{S}}}}{\eta} 
\sqrt{\frac{2 \log(p-s)}{n}} \\ \nonumber
& \geq & 
\frac{4 c_0 \sigma (32 K^2(s, s, 3, \Sigma) \lambda_{\init} \sqrt{s})}{\eta}
\sqrt{\frac{2\log(p-s)}{n}}\\ \nonumber
& = & 
\frac{4 c_0 \sigma \tilde{\delta}_{\Sc}}{\eta}
\sqrt{\frac{2\log(p-s)}{n}}
\end{eqnarray*}
and
\begin{eqnarray*}
\lambda_n
& \leq &
\frac{16 M \sqrt{2} K(s, s, 3, \Sigma)K(s, s, 3, X) c_0 \sigma \lambda_{\init} 
\sqrt{\size{\bar{S}}}}{K(s, s, 3, X)}
 \sqrt{\frac{2 \log(p-s)}{n}} \\
& \leq & 
M c_0 \sigma 32 K^2(s, 3, 3, \Sigma) \lambda_{\init}\sqrt{s}
\sqrt{\frac{2 \log(p-s)}{n}} \\
& = &
M c_0 \sigma \tilde{\delta}_{\Sc}
\sqrt{\frac{2 \log(p-s)}{n}},
\end{eqnarray*}
where we used the fact that $K(s, 3, X) \leq K(s, s, 3, X)$.
Hence~\eqref{eq::penalty-choice} is satisfied.
In addition, for $K = \sqrt{2} K(s, s, 3, \Sigma)$, 
the sparsity condition~\eqref{eq::sparsity-condition-first} holds by 
Corollary~\ref{coro:r_n-bounds}.
Condition \eqref{eq::tragic} implies 
that the condition~\eqref{eq::tragic-0} for $\beta_{\min}$ holds, 
given~\eqref{eq::norm-S-random} and~\eqref{eq::norm-Sc-random}
and Proposition~\ref{prop:eigen-K-bound-random}.
We can then invoke Theorem~\ref{thm:first} to finish the proof
with $c_0 = \sqrt{3/2}$.
\qed

\def\Z{X}
\section{Proof of the sign recovery Lemma}
\label{sec:recovery-proof}

\subsection{Preliminaries}
We first state necessary and sufficient conditions for the event
$\sign(\hat \beta) = \sign(\beta^*)$.  Note that this is
essentially equivalent to Lemmas~$2$ and~$3$ in \cite{Wai08}.
First, for $\hat{\Sigma} = X^T X/n$, let $\scov_{RT} = \inv{n} X^T_R X_T$
be the submatrix of $\scov$ with   
rows and columns indexed by $R$ and $T$ respectively. 

\begin{lemma}
\label{lemma:KKT}
Let $\vec{b} := (\sign(\beta^*_j) w_j)_{j \in S}$.
Let $\vec{w} = (w_1, w_2, \ldots, w_p),$ where $w_j > 0, \forall j$, 
be a positive weight vector.
Assume that the matrix $X_S^T X_S$ is invertible. Then for any given
$\lambda_n > 0$ and noise vector $\e \in \R^n$, 
there exists a solution $\hat{\beta}$ for the weighted Lasso such that
$$\sign(\hat \beta) = \sign(\beta^*),$$ 
if and only if the following two conditions hold:
\begin{subeqnarray}
\label{eq:lemma-Sc}
\abs{\scov_{\Sc S} (\scov_{S, S})^{-1} \left[
\frac{\Z_S^T \e}{n} - \lambda_n \vec{b}\right] - 
\frac{\Z_{\Sc}^T \e}{n}} & \leq & \lambda_n \vec{w}_{\Sc}, \\
\label{eq:lemma-S}
\sign\left(\beta^*_S + (\scov_{S S})^{-1} 
\left[\frac{\Z_S^T \e}{n} - \lambda_n \vec{b}\right]
\right) & = & \sign(\beta^*_S). 
\end{subeqnarray}
Finally, if~\eqref{eq:lemma-Sc} holds with strict inequality, then 
the solution of the weighted Lasso is unique.
\end{lemma}

\begin{proof}
Recall that we observe $Y = X \beta^* + \e$ and
$\vec{b} := (\sign(\beta^*_i) w_i)_{i \in S}$.
Let $w = (w_1, w_2, \ldots, w_p)$ be the weight vector.

First observe that the KKT conditions imply that 
$\hat{\beta} \in \R^p$ is a solution, if and only if there exists a subgradient
$$\vec{g}  \in \partial \sum_{j=1}^p w_j |\hat{\beta}_j|
 = \{z \in \R^p | \; z_i = \sign(\hat{\beta}) w_i 
\text{ for } \hat{\beta}_i \not= 0, \text{ and } |z_j| \leq w_j 
\text{ otherwise} \}$$
such that
\beq
\inv{n} X^T X \hat{\beta} - \inv{n} X^T Y + \lambda_n \vec{g} = 0,
\eeq
which is equivalent to the following linear system
by substituting $Y = X \beta^* + \e$ and re-arranging:
\begin{gather}
\label{eq:opt-beta}
\scov (\hat{\beta} - \beta^*) - \inv{n} X^T \e + 
\lambda_n \vec{g} = 0.
\end{gather}
Hence, given $X, \beta^*, \e$ and $\lambda_n >0$ the event
$\sign(\hat\beta) = \sign(\beta^*_S)$ holds
if and only if
\begin{enumerate}
\item
there exist a point $\hat{\beta} \in \R^p$ and a subgradient 
 $\vec{g} \in \partial  \sum_{j=1}^p w_j |\hat{\beta}_j|$ such that 
~(\ref{eq:opt-beta}) holds, and
\item
$\sign(\hat{\beta}_S) = \sign(\beta^*_S)$ and 
$\hat{\beta}_{\Sc} = \beta^*_{\Sc} = 0$, which 
implies that
$\vec{g}_{S} = \vec{b}$ and
$\abs{\vec{g}_{j}} \leq w_j \forall j \in \Sc$
by definition of $\vec{g}$.
\end{enumerate}

Plugging $\hat{\beta}_{\Sc} = \beta^*_{\Sc} = 0$ and
$\vec{g}_{S} = \vec{b}$ in~(\ref{eq:opt-beta}) shows that
$\sign(\hat\beta) = \sign(\beta^*)$ if and only if
\begin{enumerate}
\item
there exists a point $\hat{\beta} \in \R^p$
and a subgradient $\vec{g} \in \partial \sum_{j=1}^p w_j |\hat{\beta}_j|$
such that
\begin{subeqnarray}
\label{eq:Sc}
\scov_{\Sc S} (\hat{\beta}_S - \beta_S^*) - 
\frac{X_{S^c}^T \e}{n}&  = &  -\lambda_n \vec{g}_{S^c},\\
\label{eq:S}
\scov_{S S} (\hat{\beta}_S - \beta_S^*) - 
\frac{X_{S}^T \e}{n} =  -\lambda_n \vec{g}_{S} 
& = & -\lambda_n \vec{b},
\end{subeqnarray}
\item
and $\sign(\hat{\beta}_S) = \sign(\beta^*_S)$ and
$\hat{\beta}_{S^c} = \beta^*_{S^c} = 0$.
\end{enumerate}
Using invertibility of $X_S^T X_S$, we can solve 
for $\hat{\beta}_{S}$ and  $\vec{g}_{S^c}$ using~(\ref{eq:Sc}) and
~(\ref{eq:S}) to obtain
\begin{eqnarray*}
- \lambda_n \vec{g}_{S^c} & = &
\scov_{S^cS} (\scov_{S S})^{-1} \left[
\frac{X_S^T \e}{n} - \lambda_n \vec{b}\right] 
- \frac{X_{S^c}^T \e}{n}, \\
\hat{\beta}_S & = & \beta^*_S + (\scov_{S S})^{-1} 
\left[\inv{n}X_S^T \e - \lambda_n \vec{b}\right].
\end{eqnarray*}
Thus, given invertibility of $X_S^T X_S$,
$\sign(\hat\beta) = \sign(\beta^*)$ holds if and only if
\begin{enumerate}
\item
there exists simultaneously a point $\hat{\beta} \in \R^p$ and a 
subgradient $\vec{g} \in \partial \sum_{j=1}^p w_j |\hat{\beta}_j|$
such that
\begin{subeqnarray}
\label{eq:last-set-a}
- \lambda_n \vec{g}_{S^c}  & = &
\scov_{S^cS} (\scov_{S S})^{-1} 
\left[
\frac{X_S^T \e}{n} - \lambda_n \vec{b}\right] - \frac{X_{S^c}^T \e}{n}, \\
\label{eq:last-set-b}
\hat{\beta}_S & = & \beta^*_S +  (\scov_{S S})^{-1}  
\left[\frac{X_S^T \e}{n} - \lambda_n \vec{b}\right],
\end{subeqnarray}
\item
and $\sign(\hat{\beta}_S) = \sign(\beta^*_S)$ and 
$\hat{\beta}_{S^c} = \beta^*_{S^c} = 0$.
\end{enumerate}
Thus, for $\sign(\hat\beta) = \sign(\beta^*)$ to hold, there exists
simultaneously a point $\hat{\beta} \in \R^p$  
and a subgradient $\vec{g} \in \partial \sum_{j=1}^p w_j |\hat{\beta}_j|$
such that
\begin{eqnarray*}
\abs{\scov_{S^cS} (\scov_{S S})^{-1} \left[
\frac{X_S^T \e}{n} - \lambda_n \vec{b}\right] - 
\frac{X_{S^c}^T \e}{n}} & = &
\abs{-\lambda_n \vec{g}_{S^c}} \leq \lambda_n  \vec{w}_{S^c},\\
\sign(\hat{\beta}_S)  = 
\sign\left(\beta^*_S + (\scov_{S S})^{-1} 
\left[\inv{n}X_S^T \e - \lambda_n \vec{b}\right]\right) 
& = & \sign(\beta^*_S),
\end{eqnarray*}
given that $\abs{\vec{g}_{S^c}} \leq \vec{w}_{S^c}$ by definition of $\vec{g}$.
Thus~(\ref{eq:lemma-Sc}) and~(\ref{eq:lemma-S}) hold for the given
$X, \beta^*, \e$ and $\lambda_n >0$.
Thus we have shown the lemma in one direction.

For the reverse direction, given $X, \beta^*, \e$, 
and suppose that~(\ref{eq:lemma-Sc}) and~(\ref{eq:lemma-S}) hold for 
some $\lambda_n > 0$, we first construct a point 
$\hat{\beta} \in \R^p$ by letting
$\hat{\beta}_{S^c} = \beta^*_{S^c} = 0$ and
\begin{eqnarray*}
\hat{\beta}_{S} & = & \beta^*_S +(\scov_{S S})^{-1} 
\left[\inv{n}X_S^T \e - \lambda_n \vec{b}\right]
\end{eqnarray*}
which guarantees that
$$\sign(\hat{\beta}_S) 
 =  \sign\left(\beta^*_S + (\scov_{S S})^{-1} 
\left[\inv{n}X_S^T \e - \lambda_n \vec{b}\right]\right) 
 =  \sign(\beta^*_S)$$
by~(\ref{eq:lemma-S}).
We simultaneously construct $\vec{g}$ by letting
$\vec{g}_{S} = \vec{b}$ and
\begin{gather}
\vec{g}_{S^c} = - \inv{\lambda_n}
\left(\scov_{\Sc S} (\scov_{S S})^{-1} \left[
\inv{n}X_S^T \e - \lambda_n \vec{b}\right] - 
\inv{n} X_{S^c}^T \e \right),
\end{gather}
which guarantees that $\abs{\vec{g}_{j}} \leq w_j, \forall j \in S^c$ 
due to~(\ref{eq:lemma-S});
hence $\vec{g} \in \partial \sum_{j=1}^p w_j |\hat{\beta}_j|$.
Thus, we have found a point $\hat{\beta} \in \R^p$ and a 
subgradient $\vec{g} \in \partial \sum_{j=1}^p w_j |\hat{\beta}_j|$
such that
 $\sign(\hat{\beta}) = \sign(\beta^*)$ and the set of 
equations~(\ref{eq:last-set-a}) and~(\ref{eq:last-set-b}) is satisfied.
Hence, by invertibility of $X_S^T X_S$, 
$\sign(\hat\beta) = \sign(\beta^*)$ for the given $X, \beta^*, \e, \lambda_n$.
\end{proof}

\subsection{Uniqueness of solution}\label{subsec.uniqueness}
Finally, the uniqueness proof follows a similar argument in the revised
draft of~\cite{Wai08}. We omit the details. 
In fact, it is illustrative to rewrite the adaptive (or weighted) Lasso
program as  
follows: Let $W = \diag(w_1, \ldots, w_p)$, for $w_j > 0$, and
let the solution to~\eqref{eq::weighted-lasso} be
$$\hat \beta = W^{-1} \hat \beta_0, \; \; \; \text { where }$$ 
\begin{eqnarray}
\hat \beta_0 & := & 
\label{eq::weighted-lasso-final} 
\arg \min_{\beta_0} \inv{2n} \| Y- X W^{-1} \beta_0 \|^2_2 +  
\lambda_n \norm{\beta_0}_1 .
\end{eqnarray}
Now we can just take $X W^{-1}$ as the design matrix and 
$\beta_0 := W \beta$ as the sparse vector that we recover through 
$\hat{\beta_0}$,
by solving the standard Lasso problem as in~\eqref{eq::weighted-lasso-final}.
It is clear that uniqueness of $\hat\beta_0$ 
to~\eqref{eq::weighted-lasso-final} is equivalent to uniqueness of 
$\hat\beta$ as $W$ is a positive-definite matrix.

\subsection{Proof of Lemma~\ref{lemma:recovery}} \label{subsec.lemma-recovery} 
Let $e_i \in \R^s$ be the vector with $1$
in $i^{th}$ position and zero elsewhere; hence $\twonorm{e_i} = 1$. 

We first define a set of random variables that are relevant for
~(\ref{eq:lemma-Sc}) and~(\ref{eq:lemma-S}):
\begin{eqnarray*}
\forall j \in \Sc, \hspace{2mm}
\label{eq::V-define}
V_j & := &
\Z_j^T \Z_S (\Z_S^T \Z_S)^{-1} \lambda_n \vec{b} + 
\Z_j^T \left\{I_{n \times n} - \Z_S (\Z_S^T \Z_S)^{-1} \Z_S^T \right\}\frac{\e}{n}, \\
\label{eq::U-define}
\forall i \in S,
\hspace{2mm}
U_i & := & e_i^T
\left(\onen \Z_S^T \Z_S\right)^{-1} 
\left[\onen \Z_S^T \e - \lambda_n \vec{b}\right].
\end{eqnarray*}
Condition~(\ref{eq:lemma-Sc}) holds if and only if the event
\begin{equation*}
\label{eq::event-V}
\event(V) := \left\{\forall j \in \Sc, \; \abs{V_j} \leq 
\lambda_n w_j \right\}
\end{equation*}
is true. For Condition~(\ref{eq:lemma-S}), the event
\begin{equation*}
\label{eq::event-U}
\event(U) := \left\{\max_{i \in S} \abs{U_i} \leq \beta_{\min} \right\},
\end{equation*}
is sufficient to guarantee that
Condition~(\ref{eq:lemma-S}) holds.

We first prove that $\prob{\event(V)}$ and $\prob{\event(U)}$ both are
large. 

{\bf{Analysis of $\event(V)$}.}\enspace
Note that
$$\mu_j  = \expct{V_j} = \lambda_n  \Z_j^T \Z_S (\Z_S^T \Z_S)^{-1}
\vec{b},\ \ j \in S^c.$$
By \eqref{eq:eta}, we have $\forall j \in S^c$,
\beq
|\mu_j| \leq \lambda_n w_j (1 - \eta).
\eeq
Denote by $P = \Z_S (\Z_S^T \Z_S)^{-1} \Z_S^T = P^2$ the projection
matrix. Let 
\beq
\tilde{V_j} = \Z_j^T \left\{
\left[I_{n \times n} - \Z_S (\Z_S^T \Z_S)^{-1} \Z_S^T \right] \frac{\e}{n}
\right\},\ \ j \in S^c
\eeq
which is a zero-mean Gaussian random variable with variance
\begin{eqnarray*}
\label{eq:var-V-phi} 
\var(\tilde{V_j}) =
\frac{\sigma^2}{n^2}X_j^T\left\{
\left[\left(I_{n \times n} - P \right) \right]
\left[\left(I_{n \times n} - P \right) \right]^T \right\}X_j 
\leq \frac{\sigma^2}{n^2} \twonorm{X_j}^2
= \frac{\sigma^2 c_0^2}{n}
\end{eqnarray*}
since $\twonorm{I - P} \leq 1$.
Using the tail bound for a Gaussian random variable
\begin{eqnarray}
\label{eq::Gaus-tail-bound}
\prob{\abs{\tilde{V}_j} \geq t}
& \leq & \frac{\sqrt{\var(\tilde{V_j})}}{t}
\exp\left(\frac{-t^2}{2 \var(\tilde{V_j})}\right) \\ \nonumber
& \leq &
\frac{\sigma c_0}{\sqrt{n}t}
\exp\left(\frac{-n t^2}{2 \sigma^2 c_0^2}\right),
\end{eqnarray}
with $t = \frac{\eta \lambda_n w_{\min}(\Sc)}{2} 
\geq 2 c_0 \sigma \sqrt{\frac{2 \log(p-s)}{n}}$
and the union bound, we have
\ben
\prob{\max_{j\in S^c} \abs{\tilde{V}_j} \geq \frac{\eta \lambda_n w_{\min}(\Sc)}{2}}
& \leq & 
\frac{(p-s)\exp\left(- 4 \log (p-s)\right)}{2 \sqrt{2 \log(p-s)}} \\
& \leq & 
\inv{2 (p-s)^3 \sqrt{2 \log(p-s)}} \ .
\een
Thus, with probability at least $1- \inv{2 (p-s)^3}$,
\begin{eqnarray*}
\forall j \in \Sc, \hspace{0.5cm}
|V_j| & \leq & |\mu_j| + |\tilde{V}_j| 
\leq 
\lambda_n w_j (1 - \eta) + \frac{\eta \lambda_n w_{\min}(\Sc)}{2} \\
&\leq & \lambda_n w_j (1 - \eta/2),
\end{eqnarray*}
and $\event(V)$ holds; in fact, it holds with straight 
inequality for $\eta > 0$.

{\bf{Analysis of $\event(U)$}.}\enspace
By the triangle inequality, and on the set $\T$, 
\ben
\max_{i \in S} \abs{U_i} & \leq & 
\norm{\left(\Z_S^T \Z_S /n\right)^{-1}}_{\infty}
\norm{X_S^T \e /n}_{\infty} + 
\norm{\left( \Z_S^T \Z_S /n\right)^{-1}}_{\infty}
\lambda_n w_{\max} \\
& \leq & 
\frac{\sqrt{s}}{\Lambda_{\min}(s)} 
\left(c_0 \sigma\sqrt{24\log p /n} + \lambda_n w_{\max}(S)\right)
< \beta_{\min}, 
\een
where
$$\norm{\left({X_S^T X_S}/{n}\right)^{-1}}_{\infty} 
\leq \sqrt{s} \norm{ \left({X_S^T X_S}/{n}\right)^{-1}}_2
= \frac{\sqrt{s}}{\Lambda_{\min} \left({X_S^T X_S}/{n}\right)}
\leq  \frac{\sqrt{s}}{\Lambda_{\min}(s)},$$
by standard matrix norm comparison results and the 
restricted eigenvalue assumption. Hence, $\event(U)$ holds on the set
$\T$. 
Denote by $\F = \event(U)^c \cup \event(V)^c$. Then we have
\ben
\prob{\F} & = &  \prob{\F \cap \T^c} + \prob{\F \cap \T} \\
& \leq & \prob{\T^c} + \prob{\event(V)^c \cap \T} \\
& \leq & \prob{\T^c} + \prob{\event(V)^c} \leq  2/p^2
\een
by Lemma~\ref{lemma:gaussian-noise} and the analysis of $\event(U)$
and as $\event(V)$ above.

\section{Proof of Theorem~\ref{thm:first}}\label{finalproof}
We note that for a fixed design $X$, once we finish checking that the 
incoherence conditions and conditions on $\lambda_n$ and $\beta_{\min}$ as 
in~\eqref{eq::beta-min-recovery} are satisfied, we can then invoke 
Lemma~\ref{lemma:recovery} to finish the theorem.
For a random design, our proof follows the case of a fixed 
design after we exclude the bad event $\X^c$ for $\X$ as defined
in~\eqref{eq::good-random-design}. 
We now show that on $\X \cap \T$, where for a fixed design 
$\X^c = \emptyset$, all conditions in Lemma~\ref{lemma:recovery} 
for $c_0^2 = 3/2$ are indeed satisfied. 

First by Lemma~\ref{lemma:norm-bounds-random}, we have 
$\Lambda_{\min}(X_S^T X_S/n) \geq \Lambda_{\min}(s)$ and 
hence~\eqref{eq:incob} hold under $\X \cap \T$, 
given~\eqref{eq::eigen-cond-random}.
Now we have by $\beta_{\min} \geq 2 \tilde{\delta}_S \geq 
2 \norm{\delta_S}_{\infty}$,
\begin{eqnarray*}
\label{eq::S-init-wights-0}
\forall \ j \in S, \; \;
|\beta_{j, \init}| & \geq & \beta_{\min} - \norm{\delta_S}_{\infty} 
\geq \frac{\beta_{\min}}{2} \; \; \text{and hence } \; \; \\
\label{eq::w-max-0}
w_{\max} & \leq & \max\left\{\frac{2}{\beta_{\min}}, 1 \right\}.
\end{eqnarray*}
It also holds by 
$1 > \tilde{\delta}_{\Sc} \geq \norm{\delta_{\Sc}}_{\infty}$ 
\ben
\forall j \in \Sc, \; \; 
|\beta_{j, {\init}} | & \leq &  \norm{\delta_{\Sc}}_{\infty} 
\leq \tilde{\delta}_{\Sc} < 1
\; \; \text { and } \; \;
w_{\min} \geq  \inv{\norm{\delta_{\Sc}}_{\infty}}.
\een
Hence the choice of $\lambda_n$ in~\eqref{eq::penalty-choice} 
guarantees that
\begin{eqnarray*}
\lambda_n w_{\min} 
\geq 
\frac{\lambda_n} {\norm{\delta_{\Sc}}_{\infty}}
\geq
\frac{\lambda_n} {\tilde{\delta}_{\Sc}}
\geq 
\frac{4 c_0 \sigma}{\eta}\sqrt{\frac{2\log(p-s)}{n}}.
\end{eqnarray*}
We now show that the incoherence condition as 
in~\eqref{eq:eta-stronger} holds given $\tilde{r_n} \geq r_n$.
\begin{enumerate}
\item
Suppose $\beta_{\min} \leq 2$ satisfies~\eqref{eq::tragic-0}, we have 
$w_{\max} = 2/\beta_{\min}$ and hence
\begin{eqnarray}
\label{eq::weights-ratio-0}
\frac{w_{\min}(1 - \eta)}{w_{\max}}
\geq 
\frac{\beta_{\min} (1 - \eta)}{2\norm{\delta_{\Sc}}_{\infty}}
\geq
\frac{\beta_{\min} (1 - \eta)}{2\tilde{\delta}_{\Sc}}
\geq \tilde{r_n} \geq r_n.
\end{eqnarray}
\item
Suppose $\beta_{\min} > 2$: then $w_{\max}(S) = 1$ and 
by assumption,
\ben
\label{eq::weights-ratio-remark}
\frac{w_{\min} (1-\eta)}{w_{\max}} 
\geq \frac{1-\eta}{\norm{\delta_{\Sc}}_{\infty}}
\geq \frac{1-\eta}{\tilde{\delta}_{\Sc}} \geq \tilde{r_n}
\geq r_n.
\een
\end{enumerate}
It is clear that \eqref{eq::beta-min-recovery} is satisfied 
given~\eqref{eq::tragic-0}, if
\begin{eqnarray}
\label{eq::con-beta-0}
\beta_{\min} & \geq & 
\max\left\{\frac{4 \lambda_n \sqrt{s}}{\beta_{\min} \Lambda_{\min}(s)},
\frac{2 \lambda_n \sqrt{s}}{\Lambda_{\min}(s)}\right\}.
\end{eqnarray}
We only need to be concerned with the first term:
given the last two terms in the $\beta_{\min}$ bound, we have
\begin{eqnarray*}
\label{eq::con-beta-1}
\beta_{\min}^2 & \geq & 
\frac{4 M c_0 \sigma \tilde{\delta}_{\Sc} 
\sqrt{s}}{\Lambda_{\min}(s)}\sqrt{\frac{2 \log p}{n}}
  \; \; \; \text{ hence }\\ \nonumber
\beta_{\min} & \geq & 
\frac{4 \sqrt{s}}{\beta_{\min} \Lambda_{\min}(s)}
M c_0 \sigma  \tilde{\delta}_{\Sc}
\sqrt{\frac{2\log(p-s)}{n}}
\geq \frac{4 \lambda_{n} \sqrt{s}}{\beta_{\min} \Lambda_{\min}(s)}.
\end{eqnarray*}
Finally, we have for both fixed and random designs, let 
$\F$ be a shorthand for the event 
$\sign(\hat\beta) \not= \sign(\beta^*)$. We have
\ben
\prob{\F} \leq  \prob{(\T \cap \X)^c} + \prob{\F \cap \T \cap \X} \leq  \prob{(\T
  \cap \X)^c} + 1/p^2,
\een
where $\X^c = \emptyset$ for a fixed design, and
the last term has been bounded using Lemma~\ref{lemma:recovery} for a fixed
design or conditioned on a random design on the set $\X$ with $c_0 =
\sqrt{3/2}$. 

\bibliographystyle{ims}

\bibliography{./local}

\end{document}

%% file: macros.tex
\newcommand{\R}{\ensuremath{\mathbb R}}

\newcommand{\Z}{\ensuremath{\mathbb Z}}

\newcommand{\F}{\ensuremath{\mathcal F}}

\newcommand{\Sc}{\ensuremath{S^c}}

\newcommand{\prob}[1]{\ensuremath{{\mathbb P}\left(#1\right)}}

\newcommand{\expct}[1]{\ensuremath{{\mathbb E}\left(#1\right)}}

\newcommand{\size}[1]{\ensuremath{\left|#1\right|}}

\newcommand{\argmin}{\operatorname{argmin}}

\newcommand{\e}{\epsilon}

\newcommand{\half}{\ensuremath{\frac{1}{2}}}
\newcommand{\onen}{\textstyle \frac{1}{n}}

\newcommand{\silent}[1]{}

\newcommand{\ip}[1]{\;\langle{\,#1\,}\rangle\;}

\newcommand{\event}{{\mathcal E}}

\newcommand{\T}{{\mathcal T}}

\newcommand{\X}{{\mathcal X}}

\newcommand{\var}{\textsf{Var}}

\newcommand{\sign}{\text{\rm sgn}}
\newcommand{\init}{\text{\rm init}}

\newcommand{\inv}[1]{\frac{1}{#1}}
\newcommand{\abs}[1]{\left\lvert#1\right\rvert}
\newcommand{\twonorm}[1]{\left\lVert#1\right\rVert_2}
\newcommand{\norm}[1]{\left\lVert#1\right\rVert}

\newcommand{\func}[1]{\ensuremath{\mathrm{#1}}}
\newcommand{\diag}{\func{diag}}

\newcommand{\scov}{\hat{\Sigma}}
\newcommand{\ben}{\begin{eqnarray*}}
\newcommand{\een}{\end{eqnarray*}}

\newcommand{\beq}{\begin{equation}}
\newcommand{\eeq}{\end{equation}}